\documentclass[aos,MSNbibl,dvips]{arximspdf}
\usepackage{graphicx}
\doi{10.1214/13-AOS1178} %kopijuoti is PTS
\volume{41}
\issue{6}
\pubyear{2013}
\firstpage{3074}
\lastpage{3110}

\makeatletter

\newcommand{\rright}{\right}
\newcommand{\lleft}{\left}
\newtheorem{theorem}{Theorem}%[section]
\newtheorem{proposition}{Proposition}
\newtheorem{lemma}{Lemma}
\newproclaim{assumption}{Assumption}
\newproclaim{remark}{Remark}
\newcommand{\indep}{\perp\!\!\!\perp}
\newcommand{\whV}{\widehat{\mathbf{V}}}
\newcommand{\whT}{\widehat{\bolds{\Theta}}}
\makeatother

\begin{document}
\begin{frontmatter}

\title{Sparse PCA: Optimal rates and adaptive estimation}
\runtitle{Sparse PCA}

\begin{aug}
\author[A]{\fnms{T. Tony} \snm{Cai}\ead[label=e1]{tcai@wharton.upenn.edu}\ead[label=u1,url]{http://www-stat.wharton.upenn.edu/\textasciitilde tcai}\thanksref{t1}},
\author[A]{\fnms{Zongming} \snm{Ma}\corref{}\ead[label=e2]{zongming@wharton.upenn.edu}\ead[label=u2,url]{http://www-stat.wharton.upenn.edu/\textasciitilde zongming}\thanksref{t2}}
\and
\author[B]{\fnms{Yihong} \snm{Wu}\thanksref{t3}\ead[label=e3]{yihongwu@illinois.edu}\ead[label=u3,url]{http://www.ifp.illinois.edu/\textasciitilde yihongwu}}
\runauthor{T. T. Cai, Z. Ma and Y. Wu} \affiliation{University of
Pennsylvania, University of Pennsylvania and\\ University of Illinois
at Urbana-Champaign}
\address[A]{T. T. Cai\\
Z. Ma\\
Department of Statistics\\
The Wharton School\\
University of Pennsylvania\\
Philadelphia, Pennsylvania 19104\\
USA\\
\printead{e1}\\
\phantom{E-mail:\ }\printead*{e2}\\
\printead{u1}\\
\printead*{u2}} %adresu isvedimo komanda gale!
\address[B]{Y. Wu\\
Department of Electrical\\
\quad and Computer Engineering\\
University of Illinois\\
\quad  Urbana-Champaign\\
Urbana, Illinois 61801\\
USA\\
\printead{e3}\\
\printead{u3}}
\end{aug}
\thankstext{t1}{Supported in part by NSF FRG Grant DMS-08-54973, NSF
Grant DMS-12-08982 and NIH Grant R01 CA 127334-05.}
\thankstext{t2}{Supported in part by the Dean's Research Fund of the
Wharton School.}
\thankstext{t3}{Supported in part by NSF FRG Grant DMS-08-54973 when he
was a postdoctoral fellow at the University of Pennsylvania.}

% HISTORY:
\received{\smonth{11} \syear{2012}}
\revised{\smonth{5} \syear{2013}}

% ABSTRACT
%
\begin{abstract}
Principal component analysis (PCA) is one of the most commonly used
statistical procedures with a wide range of applications. This paper
considers both minimax and adaptive estimation of the principal
subspace in the high dimensional setting. Under mild technical
conditions, we first establish the optimal rates of convergence for
estimating the principal subspace which are sharp with respect to all
the parameters, thus providing a complete characterization of the
difficulty of the estimation problem in term of the convergence rate.
The lower bound is obtained by calculating the local metric entropy and
an application of Fano's lemma. The rate optimal estimator is
constructed using aggregation, which, however, might not be
computationally feasible.

We then introduce an adaptive procedure for estimating the principal
subspace which is fully data driven and can be computed efficiently. It
is shown that the estimator attains the optimal rates of convergence
simultaneously over a large collection of the parameter spaces. A key
idea in our construction is a reduction scheme which reduces the sparse
PCA problem to a high-dimensional multivariate regression problem. This
method is potentially also useful for other related problems.
\end{abstract}

% KEYWORDS

\end{frontmatter}

\setcounter{footnote}{3}
%s1 #&#
\section{Introduction}
\label{secintro} Due to dramatic advances in science and technology,
high-dimensional data are now routinely collected in a wide range of
fields including genomics, signal processing,\vadjust{\goodbreak} risk management and
portfolio allocation. In many applications, the signal of interest lies
in a subspace of much lower dimension and the between-sample variation
is determined by a small number of factors. For example, in
spectroscopy, the variation of the infrared and ultraviolet spectra is
driven by the concentration levels of a small number of chemical
components in the system \cite{Varmuza09}. In financial econometrics,
it is commonly believed that the variation in asset returns is driven
by a small number of common factors combined with random noise
\cite{Chamberlain83}.

Principal component analysis (PCA) is one of the most commonly used
techniques in multivariate analysis for dimension reduction and feature
extraction, and is particularly well suited for the settings where the
data is high-dimensional but the signal has a low-dimensional
structure. PCA has a wide array of applications, ranging from image
recognition to data compression to clustering. In the conventional
setting where the dimension of the data is relatively small compared
with the sample size, the principal eigenvectors of the covariance
matrix is typically estimated by the leading eigenvectors of the sample
covariance matrix which are consistent when the dimension $p$ is fixed,
and the sample size $n$ increases \cite{Anderson03}. However, in the
high-dimensional setting where $p$ can be much larger than $n$, this
approach leads to very poor estimates. At various levels of rigor and
generality, a~series of papers
%Reimann96,
\cite{Hoyle04,Baik06,Paul07,Nadler08,JohnstoneLu09,Jung09,Birnbaum12}
showed that the sample principal eigenvectors are no longer consistent
estimates of their population counterparts. For example, Baik and Silverstein \cite{Baik06}
and Paul \cite{Paul07} showed that if $p/n\to\gamma\in(0,1)$ as
$n\to\infty$, and the largest eigenvalue $\lambda_1 \leq\sqrt {\gamma}$
and is of unit multiplicity, then the leading sample principal
eigenvector $\hat{\mathbf{v}}_1$ is asymptotically almost surely
orthogonal to the leading population eigenvector $\mathbf{v}_1$, that
is, $|\mathbf{v}_1'\hat{\mathbf{v}}_1|\to0$ almost surely. Thus, in
this case, $\hat{\mathbf{v}}_1$ is not useful at all as an estimate of
$\mathbf{v}_1$. Even when $\lambda_1 > \sqrt{\gamma}$, the angle
between $\mathbf{v}_1$ and $\hat{\mathbf{v}}_1$ still does not converge
to zero unless $\lambda_1\to\infty$. In addition to being inconsistent,
sample principal eigenvectors have nonzero loadings in all the
coordinates. This renders their interpretation difficult when the
dimension $p$ is large.

%s1.1 #&#
\subsection{Sparse PCA}

In view of the above negative results in the high-dimensio\-nal
setting, a natural approach to principal component analysis in high
dimensions is to impose certain structural constraint on the leading
eigenvectors. One of the most popular assumptions is that the leading
eigenvectors have a certain type of sparsity. In this case, the problem
is commonly referred to as \emph{sparse PCA} in the literature. The
sparsity constraint reduces the effective number of parameters and
facilitates interpretation.

Various regularized estimators of the leading eigenvectors have been
proposed in the literature. See, for example,
\cite{Jolliffe03,Zou06,dAspremont07,ShenHuang08,Solo08,Witten09,Journee10}.
Theoretical analysis has so far mainly focused on the rank-one case,
that is, estimating the leading principal eigenvector $\mathbf{v}_1$.
In this case, Johnstone and Lu \cite{JohnstoneLu09} showed that the classical PCA
performed on a selected subset of variables with the largest sample
variances leads to a consistent estimator of $\mathbf{v}_1$ if the
ordered coefficients of $\mathbf{v}_1$ have rapid decay. Shen, Shen and
Marron \cite{Shen11} and Yuan and Zhang \cite{Yuan11} proposed other consistent estimators when
$\mathbf{v}_1$ has a bounded number of nonzero coefficients.
Vu and Lei \cite{Vu12} studied the rates of convergence of estimation under
various sparsity assumptions on $\mathbf{v}_1$, and Lounici \cite{Lounici12}
further considers the minimax rates with missing data. Amini and Wainwright \cite{Amini09}
investigated the variable selection property of the methods by
\cite{JohnstoneLu09} and \cite{dAspremont07} when $\mathbf{v}_1$ has
$k$ nonzero entries all of the same magnitude.
Berthet and Rigollet \cite{Berthet12}~considered minimax detection when $\mathbf{v}_1$ has
a~bounded number of nonzeros.

More recently, for estimating a fixed number $r\geq1$ of leading
eigenvectors as $n,p\to\infty$, Birnbaum et al. \cite{Birnbaum12} studied minimax rates
of convergence and adaptive estimation of the individual leading
eigenvectors when the ordered coefficients of each eigenvector have
rapid decay. When $r>1$ and some of the leading eigenvalues have
multiplicity great than one, the individual leading eigenvectors can be
unidentifiable. On the other hand, the principal subspace spanned by
them is always uniquely defined. Ma \cite{Ma11} proposed a new method for
estimating the principal subspace and derived rates of convergence of
the estimator under similar conditions to those in~\cite{Birnbaum12}.

%s1.2 #&#
\subsection{Estimation of principal subspace}

In this paper, we focus on the estimation of the principal subspace.
Both minimax and adaptive estimation are considered. Throughout the paper,
let $\mathbf{X}$ be an $n\times p$ data matrix generated as
%
%
%e1 #&#
\begin{equation}
\label{eqmodel} \mathbf{X}= \mathbf{U}\mathbf{D}\mathbf{V}' +
\mathbf{Z}.
\end{equation}
Here $\mathbf{U}$ is the $n\times r$ random effects matrix with i.i.d. $N(0,1)$
entries, $\mathbf{D}= \operatorname{diag}(\lambda_1^{1/2},\ldots,
\lambda_r^{1/2})$ with $\lambda_1\geq\cdots\geq\lambda_r >0$,
$\mathbf{V}$ is
$p\times r$ orthonormal and $\mathbf{Z}$ has i.i.d. $N(0,\sigma^2)$ entries
which are independent of $\mathbf{U}$. Equivalently, one can think of
$\mathbf{X}$ as
an $n\times p$ matrix with rows independently drawn from the
distribution $N(0,\bolds{\Sigma})$, where the covariance matrix
$\bolds{\Sigma}$ is
given by
%
%
%e2 #&#
\begin{equation}
\bolds{\Sigma}= \operatorname{Cov}({\mathbf{X}}_{{i} *})= \mathbf {V}
\bolds{\Lambda}\mathbf{V}' + \sigma^2
\mathbf{I}_p. \label{eqspike-model}
\end{equation}
Here $\bolds{\Lambda}=\operatorname{diag}(\lambda_1,\ldots, \lambda_r)$
and $\mathbf{V}= [\mathbf{v}_1,\ldots,\mathbf{v}_r]$ is $p\times r$
with orthonormal columns. The $r$ largest eigenvalues of
$\bolds{\Sigma}$ are $\lambda_i+\sigma ^2$, $i=1, \ldots, r$, and the
rest are all equal to $\sigma^2$. The $r$ leading eigenvectors of
$\bolds{\Sigma}$ are given by the columns of~$\mathbf {V}$. Since the
spectrum of $\bolds{\Sigma}$ has $r$ spikes, the covariance structure
(\ref{eqspike-model}) is commonly known as the \emph{spiked covariance
matrix model} \cite{Johnstone01} in the literature.

The goal of the present paper is to estimate the principal subspace
$\operatorname{span}(\mathbf{V})$ based on the observation $\mathbf
{X}$. Note that the principal subspace is uniquely identified with the
associated projection matrix $\mathbf{V}\mathbf{V}'$. In addition, any
estimator could be regarded as the subspace spanned by the columns of a
matrix $\whV$ with orthonormal columns, hence uniquely identified with
its projection matrix $\whV\whV'$. Thus, estimating
$\operatorname{span}(\mathbf{V})$ is equivalent to estimating
$\mathbf{V}\mathbf{V}'$. Let $\|\cdot\|_{\mathrm{F}}$ denote the
Frobenius norm. In this paper we consider optimal and adaptive
estimation of $\operatorname {span}(\mathbf{V})$ under the loss
function
%
%
%e3 #&#
\begin{equation}
\label{eqloss} L(\mathbf{V},\whV) = \bigl\|\mathbf {V}\mathbf{V}' - \whV
\whV'\bigr\|_{\mathrm{F}}^2,
\end{equation}
which is a commonly used metric to gauge the distance between linear
subspaces. It also coincides with twice the sum of the squared sines of
the principal angles between the respective linear span.

The difficulty of estimating $\operatorname{span}(\mathbf{V})$
depends on the joint sparsity of
the columns of $\mathbf{V}$. Let $\|{\mathbf{V}}_{{j} *}\|$ denote
the Euclidean norm of
the $j$th row of $\mathbf{V}$. Order the row norms in decreasing order as
$\|{\mathbf{V}}_{{(1)} *}\| \geq\cdots\geq\|{\mathbf{V}}_{{(p)}
*}\|$. We define the
\emph{weak $\ell_q$ radius} of $\mathbf{V}$ as
%
%e4 #&#
\begin{equation}
\|\mathbf{V}\|_{q,w} \triangleq\max_{j \in[p]} j \|{
\mathbf {V}}_{{(j)} *} \|^q \label{eqlqw}
\end{equation}
and let
%
%
%e5 #&#
\begin{equation}
O(p,r) = \bigl\{ \mathbf{V}\in\mathbb{R}^{p \times r}\dvtx\mathbf
{V}'\mathbf{V}=\mathbf{I}_r \bigr\} \label{eqOpr}
\end{equation}
denote the collection of $p \times r$ matrices with orthonormal columns.
We consider the following parameter spaces for $\bolds{\Sigma}$ where the
weak $\ell_q$ radius of $\mathbf{V}$ is constrained:
%
%
%e6 #&#
\begin{eqnarray}\label{paraspace}
\Theta_q(s,p,r,\lambda)
&=& \bigl\{\bolds{\Sigma}= \mathbf{V}\bolds{\Lambda}\mathbf{V}'
+ \mathbf{I}_p\dvtx0<\lambda\leq \lambda_r\leq\cdots\leq
\lambda_1 \leq\kappa\lambda,
\nonumber\\[-8pt]\\[-8pt]
&&\hspace*{98pt} \mathbf{V}\in O(p,r), \|\mathbf{V} \|_{q,w} \leq s \bigr\},\nonumber
\end{eqnarray}
where $q\in[0,2)$ and $\kappa>1$ is a fixed constant. Note that in the
rank-one case, our structural assumption coincides with
\cite{JohnstoneLu09}, (9), or \cite{Ma11}, (3.5). In the special case
of $q=0$, the union of the column supports of $\mathbf{V}$ is of size
at most
$s$. Weak $\ell_q$-ball is a commonly used model for sparsity. See, for
example, Abramovich et al. \cite{ABDJ06} for wavelet estimation and Cai and Zhou
\cite{CZ12} for sparse covariance matrix estimation. Group sparsity is
also useful for high-dimensional regression, see, for example, Lounici
et al. \cite{Lounici11}.

Let $q\in[0,2)$ and $s > 0$. Denote the weak-$\ell_q$ ball on
$O(p,r)$ by
%
%
%e7 #&#
\begin{equation}
\label{eqrow-spaceG} \mathcal{G}_q(s,p) = \bigl\{ \mathbf{V}\in O(p,r)\dvtx
\|\mathbf{V} \|_{q,w} \leq s \bigr\},
\end{equation}
which is the parameter space of $\mathbf{V}$.
%Taking into account of the orthogonality constraint, define $
In order for $\mathcal{G}_q(s,p)$ to be nontrivial, that is, neither
empty nor the whole $O(p,r)$, the weak-$\ell_q$ radius must satisfy
(see Section~7.1 in the supplementary material~\cite{appsm} for a
proof)
%
%
%e8 #&#
\begin{equation}
\frac{2-q}{2}r \leq s \leq r^{q/2} p^{(2-q)/2}. \label{eqradius}
\end{equation}
In particular, if $q = 0$, then we have $1 \leq r \leq s \leq p$.
Throughout the paper, we assume that (\ref{eqradius}) holds.

%s1.3 #&#
\subsection{Optimal rates of convergence}

Combining the upper and lower bound results developed in
Section~\ref{secminimax}, we establish the following minimax rates of
convergence for estimating the principal subspace $\operatorname
{span}(\mathbf{V})$ under the loss (\ref{eqloss}). We focus here on the
exact sparse case of $q=0$; the optimal rates for the general case of
$q\in(0, 2)$ are given in Section~\ref{secminimax}. For two sequences
of positive numbers $a_n$~and~$b_n$, we write $a_n\gtrsim b_n$ when
$a_n\geq cb_n$ for some absolute constant $c>0$ and $a_n\lesssim b_n$
when $b_n\gtrsim a_n$. Finally, we write $a_n \asymp b_n$ when both
$a_n\gtrsim b_n$ and $a_n\lesssim b_n$ hold.
% we write $a_{n}\asymp b_{n}$ if there exist absolute constants $0< c

%
%
%th1 #&#
\begin{theorem}\label{thmintrosparse} Suppose\vspace*{-1pt} we observe data $\mathbf
{X}$ as in
(\ref{eqmodel}). Let $\frac{\lambda}{\sigma^2}\gtrsim
\sqrt{\log{n}\over n}$, $s-r \gtrsim s\wedge\log\frac{\mathrm
{e}p}{s}$ and
$n \gtrsim s \log\frac{\mathrm{e}p}{s} \vee\log\frac{\lambda
}{\sigma^2}$.
% for some sufficient large constant $C$.
The minimax\vspace*{1pt} risk for estimating the principal subspace
$\operatorname {span}(\mathbf{V})$ under the loss (\ref{eqloss})
satisfies
%
%
%e9 #&#
\begin{equation}
\label{eqopt-rate0} \inf_{\whV} \sup_{\bolds{\Sigma
}\in\Theta_0(s,p,r,\lambda
)}
\mathsf{E} \bigl\|\whV\whV' -\mathbf{V}\mathbf{V}'\bigr\|
_{\mathrm{F}}^2 \asymp \frac{\lambda/\sigma^2 + 1}{n (\lambda/\sigma^2)^2} \biggl( r(s-r) + s\log
\frac{\mathrm{e}p}{s} \biggr)
\end{equation}
as long as the right-hand side of (\ref{eqopt-rate0}) does not exceed
some absolute constant. Otherwise, there exists no consistent
estimator.
\end{theorem}
%
% It is interesting to note that the optimal rate
%through $r(s-r)$, which is the dimension of %the
% % Grassmannian manifold $G(s,r)$.
% $O(s,r)$.
% Therefore the dependence on $r$ is \emph{not} monotonic, with the
%worst case occurring at $r = {s}/{2}$.
The rate of convergence in (\ref{eqopt-rate0}) depends optimally on all
the parameters $s, p, r, n$ and $\lambda$. The result thus provides a
precise characterization of the difficulty of the principal subspace
estimation problem in terms of the minimax rates over a wide range of
parameter values.

A key step in establishing the optimal rates of convergence is the
derivation of rate-sharp minimax lower bounds. It is highly nontrivial
to obtain a lower bound which depends optimally on all parameters, in
particular the eigenvalues and the rank. Our main technical tool for
the lower bounds is based on local metric entropy
\cite{LeCam73,Birge83,YB99}, instead of the usual methods based on
explicit constructions of packing sets together with Fano's lemma used,
for example, in \cite{Paul07,Birnbaum12,Vu12}. Although the method is
abstract in nature, the advantage is that it only relies on the
analytical behavior of the metric entropy of the parameter space, thus
allowing us to sidestep constructing an explicit packing, which can be
a challenging task due to the need of fulfilling both the orthogonality
and the weak-$\ell_q$ ball constraints.

We then construct an explicit estimator using an aggregation scheme,
which is shown to attain the same rates of convergence as those of the
minimax lower bounds. The matching lower and upper bounds together
establish the optimal rates of convergence. This aggregation method can
potentially be useful for other high-dimensional sparse PCA problems as
well. Aggregation methods have been widely used and well studied in
statistics literature. See, for example, Juditsky and Nemirovski
\cite{JN00}, Yang \cite{Yang00}, Nemirovski \cite{Nemirovski00} and
Rigollet and Tsybakov \cite{RT11}. To the best of our knowledge, this
is the first application of the aggregation approach to sparse PCA
which yields optimality results.

%s1.4 #&#
\subsection{Adaptive estimation}

The rate-optimal aggregation estimator depends on the model parameters
that are usually unknown in practice and is unfortunately not
computationally feasible when $p$ is large. We then propose an adaptive
estimation procedure that is fully data driven and easily
implementable. The estimator is shown to attain the optimal rate of
convergence simultaneously over a large collection of the parameter
spaces defined in (\ref{paraspace}).

The proposed method is based on a reduction scheme. By a conditioning
argument, the original sparse PCA problem is reduced to a
high-dimensional regression problem with orthogonal design and group
sparsity on the regression coefficients. Then, we apply the model
selection penalty idea from \cite{Birge01} to construct the final
estimator.

%
%estimator in multivariate regression with orthogonal design? ...}
A key step in the reduction scheme is the construction of two new
samples in the form of (\ref{eqmodel}), which share the same
realization of the random effects $\mathbf{U}$ but have independent
copies of
the noise matrix $\mathbf{Z}$. This construction works because a common
realization of $\mathbf{U}$ is critical in guaranteeing a sufficient
signal-to-noise ratio in the resulting regression problem. In contrast,
splitting the original sample into two halves fails to achieve this
goal. On the other hand, the independence of the noise components
ensures that the regression problem has white noise structure. The
adaptivity and minimax optimality of the subspace estimator depend
heavily on those of the regression coefficient estimator. Thus, as a
byproduct of the analysis, we also show that our estimator for
regression coefficients is adaptively rate optimal under group
sparsity. To the best of our knowledge, the specific estimator and its
adaptive optimality is also new in the literature.
%s1.5 #&#
\subsection{Other related work}

The present paper is related to a fast growing literature on estimating
sparse covariance/precision matrices as well as low-rank matrices.
Significant advances have been made on optimal estimation of the whole
covariance or precision matrix. Many regularization methods, including
banding, tapering, thresholding and penalization, have been proposed.
In particular, Cai, Zhang and Zhou \cite{CZZ10}~established the optimal rate of convergence
for estimating a~class of bandable covariance matrices under the
spectral norm. Cai and Yuan \cite{CY12} proposed a block thresholding
\mbox{procedure} which is shown to be adaptively rate-optimal over a
wide range of collections of bandable covariance matrices. Bickel and Levina \cite{BJ08b}
introduced a thresholding procedure and obtained rates of convergence
for sparse covariance matrix estimation. Cai and Zhou \cite{CZ12} established the
minimax rates of convergence for estimating sparse covariance matrices
under a range of matrix norms including the spectral norm. Cai, Liu and
Zhou \cite{CLZ12} obtained the optimal rate of convergence for estimating the sparse
precision matrices.

Our work is also related to another active area of research, namely,
the recovery of low-rank matrices based on noisy observations.
Negahban and Wainwright \cite{Negahban11} studied (near) low-rank matrix recovery by
$M$-estimators under restricted strong convexity\vadjust{\goodbreak} based on the penalized
nuclear norm minimization over matrices. Koltchinskii, Lounici and Tsybakov \cite{Koltchinskii11}
considered estimation of low-rank matrices based on a trace regression
model which includes matrix completion as a special case. A nuclear
norm penalized estimator was proposed and a general sharp oracle
inequality was established. See also Recht, Fazel and Parrilo \cite{rankmin} and
Rohde and Tsybakov \cite{Rhode11}.

%s1.6 #&#
\subsection{Organization of the paper}

The rest of the paper is organized as follows. After introducing basic
notation, Section~\ref{secminimax} establishes the minimax rates of
convergence for estimating the principal subspace by obtaining matching
minimax lower and upper bounds. An aggregation estimator is constructed
and shown to be rate optimal. Section~\ref{secadaptive} introduces an
adaptive estimation procedure for the principal subspace which is fully
data driven and easily computable. It is shown that this estimator
attains the optimal rates of convergence simultaneously over a large
collection of parameter spaces. Connections to other related problems
are discussed in Section~\ref{secdiscussion}. The proofs of the main
results and key technical lemmas are given in Section~\ref{secproof}
and some additional technical arguments are contained in the supplementary material~\cite{appsm}.

%s2 #&#
\section{Minimax rates for principal subspace estimation}
\label{secminimax}

We establish in this section the minimax rates of convergence for
estimating the principal subspace in two steps. First, minimax lower
bounds are obtained for the estimation problem under the loss
(\ref{eqloss}). Then an aggregation estimator is introduced and is
shown to attain the same rates as given in the lower bounds, under mild
conditions on the parameters.
%when the parameters are appropriately chosen.
The matching lower and upper bounds thus establish the minimax rates of
convergence.
%Without loss of generality, we assume $\bsmu= {\bf0}$ from now on.
% Let $\S= \frac{1}{n}\sum_{i=1}^n \bfx_i \bfx_i'$ be the sample
%covariance matrix. \nr{Should be use the notation $\row{\X}{i}$ here?
%-- YW}

We begin by introducing some basic notation. Throughout the paper, for
any matrix $\mathbf{X}= (x_{ij})$ and any vector $\mathbf{u}$, denote
by $\|\mathbf{X}\|$ the spectral norm, $\|\mathbf{X}\|_{\mathrm{F}}$
the Frobenius norm and $\| \mathbf{u}\|$ the vector $\ell_2$ norm.
Moreover, the $i$th row of $\mathbf{X}$ is denoted by
${\mathbf{X}}_{{i} *}$ and the $j$th column by ${\mathbf{X}}_{* {j}}$.
Let $\operatorname{supp}(\mathbf{X}) = \{i\dvtx{\mathbf{X}}_{{i} *}
\neq0\}$ denote the row support of $\mathbf{X}$. For a positive integer
$p$, $[p]$ denotes the index set $\{1, 2,\ldots, p\}$. For two subsets
$I$ and $J$ of indices, denote by $\mathbf{X}_{IJ}$ the $|I|\times|J|$
submatrices formed by $x_{ij}$ with $(i,j) \in I \times J$. Let
${\mathbf{X}}_{{I} *} = \mathbf{X}_{I {[p]}}$ and ${\mathbf {X}}_{*
{J}} = \mathbf{X}_{{[n]} J}$. For any square matrix $\mathbf{A}=
(a_{ij})$, we let $\operatorname {Tr}(\mathbf{A}) = \sum_{i}a_{ii}$ be
its trace. Define the inner product of matrices $\mathbf{B}$ and
$\mathbf{C}$ of the same size by $\langle\mathbf {B}, \mathbf{C}
\rangle= \operatorname{Tr}(\mathbf{B}'\mathbf{C})$. For any matrix
$\mathbf{A}$, we use $\sigma_i(\mathbf{A})$ to denote its $i$th largest
singular value. When $\mathbf{A}$ is positive semi-definite, $\sigma
_i(\mathbf{A})$ is also the $i$th largest eigenvalue of $\mathbf{A}$.
Let $\operatorname{span}(\mathbf{A})$ denote the linear subspace
spanned by the columns of $\mathbf{A}$.
% denote the $|A| \times p$ submatrices formed by the rows $\row{
%the indices of the rows of $\X$ with at least one nonzero entries.
For any real numbers $a$ and $b$, set $a\vee b = \max\{a,b\}$ and
$a\wedge b = \min\{a,b\}$.
%positive numbers, we write $a_n\gtrsim b_n$ when $a_n\geq cb_n$ for
%some numeric constant $c$, and $a_n\lesssim b_n$ when $a_n\leq Cb_n$
%for some numeric constant $C$. Finally, we write $a_n \asymp b_n$ when
%both $a_n\gtrsim b_n$ and $a_n\lesssim b_n$ hold.}
For any set $A$, $|A|$ denotes its cardinality. Let $\mathbb{S}^{p-1}$
denote the unit sphere in $\mathbb{R}^p$. Let $G(k,r)$ denote the
Grassmannian manifold consisting of all $r$-dimensional linear subspace
of $\mathbb{R}^k$. Let $O(p)$ denote the collection of all $p \times p$
orthogonal matrices. Throughout the paper, we use $c$ and $C$ to denote
generic absolute positive constants, though the actual value may vary
at different occasions. For any sequences $\{a_n\}$ and $\{b_n\}$ of
positive numbers, we write $a_n\gtrsim b_n$ when $a_n\geq cb_n$ for
some absolute constant $c$, and $a_n\lesssim b_n$ when $a_n\leq C b_n$
for some absolute constant $C$. Finally, we write $a_n \asymp b_n$ when
both $a_n\gtrsim b_n$ and $a_n\lesssim b_n$ hold.

% The parameters $\lambda> 0$, $\kappa> 1$ are free.

% Now let $q \in(0,2)$. Note that $\calF_{q'}(s,p) \cap O(p,r) \subset
%particular, $q'=0$. Then we must have $s \geq r$.

%We also need to make sure \eqref{eqkqs} has a solution. Then we
%have $1 \leq r \leq s \leq k^* \leq p$.

%{\red Assumptions to be finalized}
% \item There exists a sufficiently small constant $\eps>0$, such that
%nh(\lambda) \geq\frac{k}{\eps^2}.
%Used in the oracle upper bound.

%nh(\lambda) \geq C(r + \log\frac{\eexp p}{k^*})

%nh(\lambda) \geq C k^* (r + \log\frac{\eexp p}{k^*})
%where $k=k^*_q$ in \eqref{eqkqs}.

%Theorem~\ref{thmlower-bound0} need assumption stated there.

%s2.1 #&#
\subsection{Lower bounds}
\label{secLB}

We first establish the minimax lower bounds which are instrumental in
obtaining the optimal rates of convergence. In view of the upper bounds
to be given in Section~\ref{secUB} by an aggregation procedure, these
lower bounds are minimax rate optimal under mild conditions.

Before proceeding to the precise statements, we introduce the following
notation: let
%
%
%e10 #&#
%e11 #&#
\begin{eqnarray}
\label{eqh} h(\lambda) &=& \frac{\lambda^2}{\lambda+1},
\\
\Psi(k,p,r,n,\lambda) &=& \frac{1}{n h(\lambda)} \biggl( rk + k\log
\frac{\mathrm{e}p}{k} \biggr)\label{eqrate}
\end{eqnarray}
and
%
%
%e12 #&#
\begin{equation}
\Psi_0(k,p,r,n,\lambda) = \frac{1}{n h(\lambda)}
\biggl( r(k-r) + k\log\frac{\mathrm{e}p}{k} \biggr).\label{eqrate-0}
\end{equation}
Define the \emph{effective dimension} by
%k^*_q(s,p,r,\lambda,n) = s \pth{\frac{n h(\lambda)}{r + \log p} }^{q/2}
% \label{eqkqs}
%
%
%e13 #&#
\begin{equation}
\label{eqkqs} k^*_q(s,p,r, n,\lambda) \triangleq{ \bigl
\lceil{x_q(s,p,r,n, \lambda )} \bigr\rceil},
\end{equation}
where ${ \lceil{a}  \rceil}$ denotes the smallest integer
no less than $a \in\mathbb{R}
$, and
%
%
%e14 #&#
\begin{equation}
\label{eqxqs} x_q(s,p,r,n,\lambda) \triangleq\max \biggl\{ 0 \leq x
\leq p\dvtx x \leq s \biggl( \frac{n h(\lambda)}{r + \log(\mathrm{e}p/x)} \biggr)^{q/2} \biggr\}.
\end{equation}
%
%--ZM
%No!!! We have used too much the max version!! Also where did we need
%$k_q^*$ to be greater than $r$? Then maybe there is no solution if one
%insists on $k\geq r$. I think I wrote some crap here. So ignore. --
%Yihong}
%%\min\sth{k \in[p]\dvtx k \geq r, \frac{4 q}{2-q}k(s/k)^{2/q} n h(

%k^*_q(s,p,r,{\blue n,\lambda}) = \floor{x_q(s,p,r,{\blue n,\lambda})}
%%\min\sth{k \in[p]\dvtx k \geq r, \frac{4 q}{2-q}k(s/k)^{2/q} n h(
% \label{eqkqs}
%where for any number $a$, $\lfloor a\rfloor$ is the largest integer
%smaller than $a$, and $x_q(s,p,r,\lambda,n)$ is the solution to the
%following equation
%x = s \pth{ \frac{n h(\lambda)}{r + \log\frac{\eexp p}{x}}}^{q/2}.
%%\min\sth{k \in[p]\dvtx k \geq r, \frac{4 q}{2-q}k(s/k)^{2/q} n h(
% \label{eqxqs}

%

\begin{figure}%[b]
\begin{tabular}{@{}cc@{}}

\includegraphics{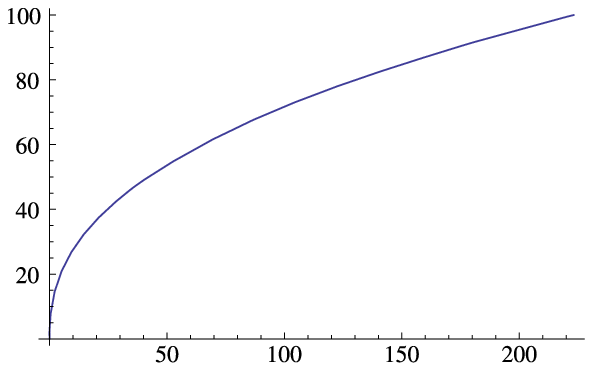}
 & \includegraphics{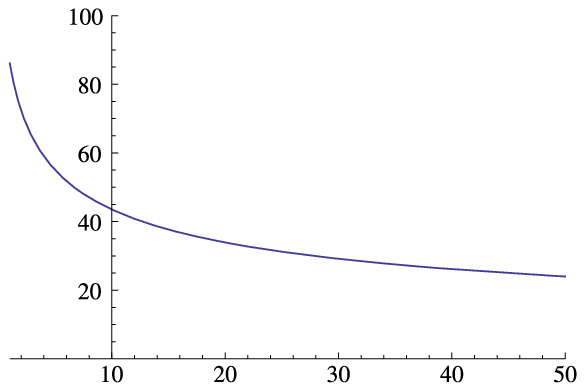}\\
\footnotesize{(a) $x_q \sim nh(\lambda)$} & \footnotesize{(b) $x_q \sim r$}\\[6pt]

\includegraphics{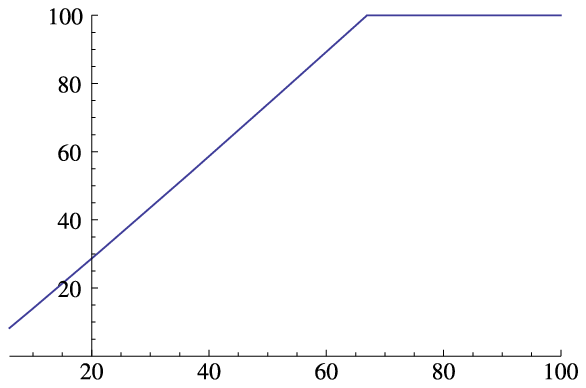}
 & \includegraphics{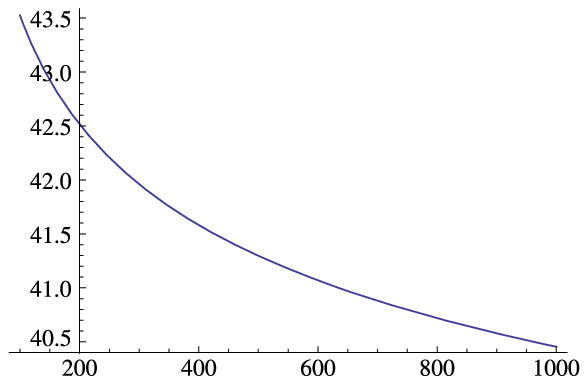}\\
\footnotesize{(c) $x_q \sim s$} & \footnotesize{(d) $x_q \sim p$}
\end{tabular}
%parameters (default values: $p = 100, s = 30, r = 10,nh(\lambda) =
%30,q = 0.8$).}
\caption{Plots of $x_q$ against individual parameters [default values:
$p = 100$, $s = 30$, $r = 10$, $nh(\lambda) = 30$, $q = 0.8$]. The effective
dimension is $k_q^*={ \lceil{x_q}  \rceil}$.}\label{figkqs}
\end{figure}

%
%
%re1 #&#
\begin{remark}[(Effective dimension)]\label{rmkeffdim}
The effective dimension $k_q^*$ is a proxy which captures the
massiveness of the parameter set for the principle subspace under the
weak-$\ell_q$ constraint.
In addition, the minimax estimation rate turns out to be a strictly
increasing function of $k_q^*$.
% Moreover, we will show that regular PCA is minimax rate optimal if
%and only if that the effective dimension is on the same order as the
%ambient dimension. See Remark~\ref{rmkregPCA} for more details.
In the exact sparse case, it is evident from~(\ref{eqkqs}) that $k^*_0
= s$. Therefore in this case, the effective dimension coincides with
the row sparsity of $\mathbf{V}$. For $q\in(0,2)$,
% In addition,
the effective dimension satisfies the following properties (proved in
Section~7.2 in the supplementary material~\cite{appsm}):
\begin{longlist}[(3)]
\item[(1)] % For any $q \in(0,2)$,
$k_q^* \geq1$.

\item[(2)] $k_q^* = p$ if and\vspace*{1pt} only if $s \geq p  ( \frac{r
    + 1}{nh(\lambda)}  )^{q/2}$, in which case the effective dimension
    coincides with the ambient dimension.

\item[(3)] $s \mapsto k_q^*$ is increasing. Moreover, there exists a
function $\tau_q$, such that $k^*_q(a s,p,r, n,\lambda) \leq
k^*_q(s,p,r, n,\lambda) \tau_q(a)$ for any $a \geq1$.

\item[(4)] $k_q^* \gtrsim s$ if and only if the assumption (\ref
{eqassumption-bigsnr}) holds. % \nr{forward reference?}.
% Under the assumption \eqref{eqassumption-bigsnr}, it can be
%shown that the solution satisfies $x_q(s,p,r, {n,\lambda}) \geq s$.
%Consequently, $k_q^*(s,p,r, n,\lambda) \geq s$.
\end{longlist}
See Figure~\ref{figkqs} for a graphical illustration on the dependence of
the effective dimension $k_q^*$ on various parameters.
% The monotonicity is evident from the definition.
%effective dimension goes down...}
%If $q \in(0,2)$, then we have $\frac{4 q}{2-q}k(s/k)^{2/q} n h(
%k_q^* \geq C s \pth{\frac{n}{r + \log p} }^{q/2}
%where $C$ only depends on $q$.
%
%and the minimax rate is
%s \pth{\frac{r + \log p}{n h(\lambda)}}^{q/2-1}.
\end{remark}

%
%f1 #&#

Without loss of generality, we assume unit noise variance ($\sigma
^2=1$) from now on. All results hold for a general $\sigma$ by
replacing $\lambda$ with $\lambda/\sigma^2$. We consider the lower
bounds separately in two cases: $0<q<2$ and $q=0$.

%
%th2 #&#
\begin{theorem}[(Lower bound: $0 < q < 2$)]\label{thmlower-bound}
Let $p \in\mathbb{N}$, $r \in[p]$ and $k^*_q$ be defined in
(\ref{eqkqs}).
% Let $s \geq1$ satisfies \eqref{eqradius}.
% \nr{The condition on $s$ is redundant in view of (15).
% agree. dude can you go ahead and revise as well as in the proof if
%any? -- yihong}
Let the observed matrix $\mathbf{X}$ be generated by model (\ref
{eqmodel}) with
$\sigma=1$.
% Let $k^*_q$ be defined in \eqref{eqkqs}.
Assume that
%
%
%e15 #&#
\begin{equation}
% k^*_q \geq2 r
r \leq\frac{s}{2} \wedge \bigl(p+1-k^*_q \bigr)
\label{eqassumption-rs2}
\end{equation}
%
%r \leq p+1-k^*_q
% \label{eqassumption-prk-q}
and that\vspace*{-3pt}
%
%
%e16 #&#
\begin{equation}
\label{eqassumption-bigsnr} %nh(\lambda) \geq C_0 \pth{r + \log\frac{
nh(\lambda) \geq C_0^{2/q}
\biggl( r + \log\frac{\mathrm{e}p}{s} \biggr)
\end{equation}
for some sufficiently large absolute constant $C_0$.
% If
% \begin{equation}
%implied by previous assumption}}
Then there exists a constant $c$~depending only on $q$ and an absolute
constant $c_0$, such that the minimax risk for estimating $\mathbf{V}$
over the parameter space $\Theta= \Theta_q(s,p,r,\lambda)$ satisfies
%
%
%e17 #&#
\begin{equation}
\inf_{\whV} \sup_{\bolds{\Sigma}\in\Theta} \mathsf{E}\bigl\|\whV
\whV' -\mathbf{V}\mathbf{V}' \bigr\|_{\mathrm{F}}^2
\geq c \Psi \bigl(k^*_q,p,r,n,\lambda \bigr) \wedge c_0,
\label{eqlower-boundq}
\end{equation}
where $\Psi$ is defined in (\ref{eqrate}).
\end{theorem}

%Note that the above lower bound is obtained under the assumption
%radius $s$, the rank of $\V$ could take values up to $\frac{2s}{2-q}$.
%In particular, in the exact sparse case where $q = 0$, $r$ takes
%values in the full range $[s]$. An intriguing question is what happens
%when the rank $r$ exceeds $\frac{s}{2}$?
%answer turns out to be interesting: In the sparse case, the
%statistical difficulty for estimating the $r$ leading singular vectors
%depend on the $r$ only through $r(s-r)$, which is the dimension of the
%Grassmannian manifold $G(s,r)$. Therefore the dependence on $r$ is
%$r$ by $s-r$. }
%The following more precise lower bound characterizes this behavior
%precisely for the case $q=0$. % Since the rank $r$ takes values in
%$[k]$,

For the case of $q=0$ we have the following lower bound:
%
%
%th3 #&#
\begin{theorem}[(Lower bound: $q = 0$)]
\label{thmlower-bound-0} Let $p,s,r$ be integers such that $1 \leq r
\leq s \leq p$. Let the observed matrix $\mathbf{X}$ be generated by model
(\ref{eqmodel}) with $\sigma=1$.
%Assume that $s$ and $r$ are positive integers satisfying
%%\begin{equation}
%%r \leq p+1-s.
%% \label{eqassumption-prk}
%%\end{equation}
% Then there exist an absolute constant $c$, such that
Then the minimax risk for estimating $\mathbf{V}$ over the parameter space
$\Theta= \Theta_0(s,p,r,\lambda)$ satisfies\vspace*{-4pt}
%
%
%e18 #&#
\begin{eqnarray}\label{eqlower-bound0}
&& \inf_{\whV} \sup_{\bolds{\Sigma}\in\Theta} \mathsf{E}\bigl\|\whV
\whV' -\mathbf{V}\mathbf{V}' \bigr\|_{\mathrm{F}}^2
\nonumber\\[-11pt]\\[-11pt]
&&\qquad \gtrsim %\frac{c}{n h(\lambda)}\pth{r(s-r) + (s-r) \log\frac{e(p-r)}{(s-r)}}
\biggl[ \frac{1}{n h(\lambda)} \biggl( r(s-r) + (s-r)
\log\frac
{\mathrm{e} (p-r)}{s-r} \biggr) \biggr] \wedge1.\nonumber
\end{eqnarray}
%If the condition \eqref{eqassumption-min-snr} is not satisfied,
%then
%c'
% \label{eqlower-bound0-trv}
%for some some absolute constant $c'$.
\end{theorem}

%s2.2 #&#
\subsection{Optimal estimation via aggregation}\label{secUB}

%Define the parameter space of $\V$ under weak-$\ell_q$ constraint by
% \label{eqGqs}
%% \{\V\dvtx\V\in O(p,r), |\supp(\V)| \leq k \}
%where $\calF(q,s)$ is defined in \eqref{eqpara-space-g} and
%$O(p,r)$. In particular,

We now show that the lower bounds given in Section~\ref{secLB} are
indeed rate optimal under mild technical conditions. The optimal
estimator of $\mathbf{V}$ is constructed using sample splitting and
aggregation. The estimator is theoretically interesting but
computationally intensive. We will construct a data-driven and easily
implementable estimator in Section~\ref{secadaptive} under stronger
conditions.

We first note that the loss function (\ref{eqloss}) satisfies
% the following
%
%
%e19 #&#
\begin{equation}
\label{eqlossp} L(\mathbf{V}, \whV) = 2 r - 2 \bigl\| \whV' \mathbf{V}
\bigr\|_{\mathrm{F}}^2 = 2 \bigl\| \bigl(\mathbf{I}- \mathbf{V}
\mathbf{V}' \bigr) \whV \whV'\bigr\|_{\mathrm{F}}^2.
\end{equation}
Moreover, the loss function is invariant under orthogonal complement,
that is, $L(\mathbf{V}, \whV) = L(\mathbf{V}^\perp, \whV^\perp)$, where
$[\mathbf{V}, \mathbf{V}^\perp], [\whV,\whV^\perp ]$ are orthogonal
matrices. Therefore the loss (\ref{eqlossp}) admits the following upper
bound: % never exceeds $2 (r\wedge(p-r))$.
%
%
%e20 #&#
\begin{equation}
L(\mathbf{V}, \whV) \leq2 \bigl(r\wedge(p-r) \bigr). \label{eqlub}
\end{equation}

For notational simplicity we assume that the sample size is $2n$ and we
split the sample equally according to $\mathbf{X}=
\bigl[{\fontsize{8.36}{10}\selectfont{\matrix{\mathbf {X}_{(1)} \cr
\mathbf{X}_{(2)}}}}\bigr]$, where $\mathbf{X}_{(i)} = \mathbf{U}_{(i)}
\mathbf{D}\mathbf{V}' + \mathbf{Z}_{(i)}, i = 1,2$.
% be the two independent samples generated by \eqref{eqsplit}
%are obtained. Two equal size subsamples?} and
Denote by $\mathbf{S}_{(i)} = \frac{1}{n} \mathbf{X}_{(i)}' \mathbf
{X}_{(i)}$ the corresponding sample covariance matrix. The main idea is
to construct a family of estimators $\{\whV_B\}$ using the first
sample, indexed by the row support $B \subset[p]$, where $\whV_B$ is
the optimal estimator one would use if one knew beforehand that
$\operatorname{supp}(\mathbf{V})=B$.
% \nr{Yihong, could you add some comment here on this ``oracle''
%optimal estimator?}
Then we aggregate these estimators by selection using the second sample.

Recall\vspace*{-2pt} the effective dimension $k^*_q$ defined in
(\ref{eqkqs}). For each $B \subset[p]$ such that $|B|=k^*_q$, we define
$\whV_B \in O(p,r)$ as the $r$ leading singular vectors of
$\mathbf{J}_B \mathbf{S}_{(1)} \mathbf{J}_B$, where $\mathbf{J}_B$ is
the diagonal matrix given by
%
%
%e21 #&#
\begin{equation}
(\mathbf{J}_B)_{ii} = {\mathbf{1}_{ \{{i\in B} \}}}.
\label{eqJB}
\end{equation}
Given the collection of the $\whV_{B}$'s, we set
%
%
%e22 #&#
\begin{equation}
B^* = \mathop{\mathop{\operatorname{argmax}}_{B \subset[p]}}_{|B| = k^*_q}
\operatorname{Tr} \bigl( \whV_B' \mathbf{S}_{(2)}
\whV_B \bigr) \label{eqBstar}
\end{equation}
and define the aggregated estimator by
%
%
%e23 #&#
\begin{equation}
\whV_* = \mathbf{V}_{B^*}. \label{eqcombi}
\end{equation}

It is natural to use the same sample covariance matrix to construct the
$\whV_B$'s and to select $B^*$. The main advantage of sample splitting
is to decouple the selection of the support and the computation of the
estimator. Thus, conditioning on the first sample, we can treat the
candidate estimators as if they are deterministic, which greatly
facilitates the analysis. Sample splitting is commonly used in
aggregation based estimation, where a sequence of estimators is
constructed from the first sample and the second sample is used to
aggregate these candidates to produce a final estimator.

Estimator (\ref{eqcombi}) requires knowledge of the value of $q$, the
weak-$\ell_q$ radius $s$, the rank $r$ and the spike size $\lambda$.
Moreover, it can be computationally intensive for large values of $p$
since in principle one needs\vspace*{-2pt} to enumerate all ${p \choose
k_q^*}$ possible column supports in order to obtain $B^*$. Nonetheless,
the next theorem establishes its minimax rate optimality:
%
%
%th4 #&#
\begin{theorem}
%For any $1 \leq r \leq p$ and $n \in\naturals$.
% Assume that \eqref{eqassumption-n-k-lambda} and
Let $q \in[0,2)$. Let $k^*_q$ be defined in (\ref{eqkqs}). Let $\whV_*$
be the aggregated estimator defined in (\ref{eqcombi}). Assume that
%
%
%e24 #&#
%e25 #&#
\begin{eqnarray}
\lambda&\geq& C_0 \sqrt{\frac{\log{n}}{n}}, \label
{eqassumption-large-lambda}
\\
\label{eqassumption-superbigsnr} nh(\lambda) &\geq& C_0
k^*_q \biggl( r + \log\frac{\mathrm{e}p}{k^*_q} \biggr)
\end{eqnarray}
and % $n \geq C_0 (k \log\frac{\eexp p}{k} \vee\log\lambda)$
%
%
%e26 #&#
\begin{equation}
n \geq C_0 \biggl( k_q^* \log\frac{\mathrm{e}p}{k_q^*} \vee\log
\lambda \biggr) \label{eqassumption-lambdan}
\end{equation}
%
% n \geq C_0 r
% \label{eqassumption-rn}
for some sufficiently large constant $C_0$.
%and there exists a sufficiently small constant $\eps>0$, such that
%nh(\lambda) \geq\frac{k^*_q}{\eps^2}.
%Used in the oracle upper bound.
Then there exists a constant $C$ depending only on $\kappa$ and $q$
such that for $\Theta= \Theta_q(s,r,p,\lambda)$,
%
%
%e27 #&#
\begin{equation}
\sup_{\bolds{\Sigma}\in\Theta} \mathsf{E}\bigl\|\whV_*\whV_*' -
\mathbf{V} \mathbf{V}'\bigr\|_{\mathrm{F}}^2 \leq2 \bigl(r
\wedge(p-r) \bigr) \wedge C\Psi \bigl(k^*_q,p,r, n,\lambda \bigr),
\label{eqagg}
\end{equation}
where $\Psi(k,p,r,n,\lambda)$ and $k^*_q$ are defined in (\ref{eqrate})
and (\ref{eqkqs}), respectively. Moreover, if $q = 0$, then $\Psi$ in
(\ref{eqagg}) can be replaced by $\Psi_0$ defined in (\ref{eqrate-0})
with $k^*_0=s$, and condition (\ref{eqassumption-superbigsnr}) can be
dropped. \label{thmagg}
\end{theorem}
%
% \nb{Double check conditions on $nh(\lambda)$!}

When $q\in(0,2)$, under the conditions of Theorems \ref{thmlower-bound}
and \ref{thmagg}, the lower and upper bounds together yield the minimax
rates of convergence $\Psi(k^*_q,p,r,n,\lambda)$ given
in~(\ref{eqrate}) with the optimal dependence on all the parameters, in
particular the eigenvalues and the rank. When $q = 0$, the lower and
upper bounds match under less restrictive conditions, which will be
discussed in more detail in Remark~\ref{rmkreqs} below.
% The results thus completely characterize the difficulty of estimating
%the principal subspace in term of the minimax rates. \nr{This
%paragraph seems out of place. Relocate to Sec. 2.3 or kill?}

%s2.3 #&#
\subsection{Comments}

We conclude this section with a few important remarks.
%Note that the above lower bound is obtained under the assumption
%radius $s$, the rank of $\V$ could take values up to $\frac{2s}{2-q}$.
%In particular, in the exact sparse case where $q = 0$, $r$ takes
%values in the full range $[s]$. An intriguing question is what happens
%when the rank $r$ exceeds $\frac{s}{2}$?
%answer turns out to be interesting: In the sparse case, the
%statistical difficulty for estimating the $r$ leading singular vectors
%depend on the $r$ only through $r(s-r)$, which is the dimension of the
%Grassmannian manifold $G(s,r)$. Therefore the dependence on $r$ is
%$r$ by $s-r$. }
%The following more precise lower bound characterizes this behavior
%precisely for the case $q=0$. % Since the rank $r$ takes values in
%$[k]$,

%
%
%re2 #&#
\begin{remark}[(Minimax rates in the exact sparse case)]\label{rmkreqs}
Comparing the lower and upper bounds for $q=0$ in Theorems
\ref{thmlower-bound} and \ref{thmagg}, we see a sufficient condition
for the minimax rate to match (and hence coincide with $\Psi_0$) is
%
%
%e28 #&#
\begin{equation}
s - r \gtrsim s \wedge\log\frac{\mathrm{e}p}{s}. \label{eqPsi0tight}
\end{equation}
To see this, suppose that $s - r \gtrsim s$. Then there exists $c \in
(0,1)$, such that $r \leq(1-c) s \leq(1-c) p$. Then $(s-r)\log\frac
{\mathrm{e}(p-r)}{s-r} \geq c s \log\frac{\mathrm{e}c p}{s} \gtrsim s
\log \frac{\mathrm{e}p}{s}$ and the lower bound in
(\ref{eqlower-bound0}) agrees with $\Psi_0$. Now\vspace*{-1pt} suppose
that $s - r \lesssim s$ and $s - r \gtrsim\log\frac{\mathrm{e}p}{s}$.
Then $r\asymp s$ and $r (s-r) + s \log\frac{\mathrm{e}p}{s} \asymp r
(s-r) \asymp s (s-r)$. Since $r
\mapsto(s-r)\log\frac{\mathrm{e}(p-r)}{s-r}$ is decreasing on $[0,s]$,
we have $r (s-r) + (s-r)\log\frac{\mathrm{e}(p-r)}{s-r} \asymp s
(s-r)$. Hence the lower bound in (\ref{eqlower-bound0}) also agrees
with $\Psi_0$.

%Note that the above lower bound is obtained under the assumption
%radius $s$, the rank of $\V$ could take values up to $\frac{2s}{2-q}$.
%In particular, in the exact sparse case where $q = 0$, $r$ takes
%values in the full range $[s]$. An intriguing question is what happens
%when the rank $r$ exceeds $\frac{s}{2}$?
It is interesting to note that under the condition (\ref{eqPsi0tight}),
the minimax rate for estimating the $r$ leading singular vectors depend
on the $r$ only through $r(s-r)$, which is the dimension of the
Grassmannian manifold $G(s,r)$. Therefore the dependence on $r$ is
\emph{not} monotonic, with the worst case happening at $r =
\frac{s}{2}$.
%Moreover, the minimax rate is invariant if we replace $r$ by $s-r$.
However, it should be noted that in order for the minimax rate to
coincide with $\Psi_0$, it is \emph{necessary} to have $r$ strictly
bounded away from $s$, for example, in the regime of
(\ref{eqPsi0tight}). When $r=s$, the lower bound in
Theorem~\ref{thmlower-bound-0} becomes zero. In this degenerate case,
the only uncertainty is in the support of $\mathbf{V}$. The minimax
rate is indeed much faster than $\Psi_0$, because in this regime the
support can be estimated accurately. See Section~7.3 in the
supplementary material~\cite{appsm}.
\end{remark}

%
%
%re3 #&#
\begin{remark}\label{rmkasymp} %[Asymptotics of minimax rates]
% Although the results in the present paper are nonasymptotics, it is
%interesting to analyze the asymptotics behavior of the minimax rates.
For $q \in(0,2)$, the minimax rate $\Psi$ depends on the effective
dimension $k_q^*$ which is defined implicitly through equations
(\ref{eqkqs})--(\ref{eqxqs}). It is possible to obtain an explicit
formula of the minimax rate in some regime. For example, if $s \geq
p^{1-\epsilon} (\frac{r + \log p}{nh(\lambda)})^{q/2}$ for some
constant $\epsilon\in(0,1)$, then the effective dimension satisfies
$k_q^* \leq p^{1-\epsilon}$. Moreover, we have $k_q^* \asymp s (\frac{n
h(\lambda)}{r + \log p })^{q/2}$. Hence the minimax rate is given by
\[
\Psi \bigl(k^*_q,p,r, n,\lambda \bigr) \asymp s \biggl(
\frac{r +
\log p }{n h(\lambda)} \biggr)^{1-q/2}.
\]
\end{remark}

An interesting side product of the proofs of Theorems
\ref{thmlower-bound-0} and \ref{thmagg} is the following nonasymptotic
minimax rate\vadjust{\goodbreak} for the regular PCA problem without structural assumptions
on the principle subspaces. It is a classical result (see, e.g.,
\cite{Stein56,Eaton70}) that when $p \leq n$, the sample covariance
matrix is not exact minimax optimal for estimating the whole covariance
matrix under certain losses (e.g., the Stein loss). As shown in the next
theorem, in the unstructured case, it turns out that the sample version
of the principle subspace is minimax \emph{rate} optimal even in high
dimensions. For more details see Theorems \ref{thmoracle-lb} and
\ref{thmoracle-ub} in Sections~\ref{secpflower-bound} and
\ref{seccombipf}.

%
%
%th5 #&#
\begin{theorem}\label{thmregularPCA}
Let $\Theta= \Theta_0(p,p,r,\lambda)$.
Let $n \geq C_0 (r + \log\lambda)$ and $\lambda\geq\break  C_0 \sqrt{{\log
(n)}/{n}}$ for some sufficiently large constant $C_0$. Then for all $r
\in[p]$,
% $n \geq1$, and $\lambda> 0$,e
%
%
%e29 #&#
\begin{equation}
\inf_{\whV}\sup_{\bolds{\Sigma}\in\Theta} \mathsf {E}\bigl\|\whV
\whV' - \mathbf{V}\mathbf {V}'\bigr\|_{\mathrm{F}}^2
\asymp r \wedge(p-r) \wedge\frac{r(p-r)}{nh(\lambda)}, \label{eqregularPCA}
\end{equation}
which can be attained by $\whV$ consisting of the $r$ leading
eigenvectors of the sample covariance matrix $\mathbf{S}$.
\end{theorem}

Theorem~\ref{thmregularPCA} implies that, without structural
assumptions on the principle subspace $\mathbf{V}$, consistent
estimators exist
if and only $\frac{nh(\lambda)}{r(p-r)} \to\infty$. Moreover, unless
$nh(\lambda)$ exceeds a constant factor of $p$,
%probably we should avoid using it.}
even the optimal estimator is within a constant factor of $r \wedge
(p-r)$, the upper bound of the loss function.

In the structured case, we can also investigate when regular PCA is
rate optimal.
%Theorem~\ref{thmregularPCA} gives the we know that the .
It is intuitive to expect that regular PCA is strictly suboptimal if
the principal eigenvectors are highly sparse, since the procedure
ignores the structure of the problem. Indeed,
Theorem~\ref{thmregularPCA}, together with Theorems
\ref{thmlower-bound}--\ref{thmagg}, reveals the precise regime where
regular PCA is minimax rate optimal: under the conditions of
Theorem~\ref{thmoracle-ub}, regular PCA achieves minimax rate if and
only if the effective dimension $k_q^*\asymp p$.
% is on the same order as the ambient dimension $p$.
In view of definition (\ref{eqkqs}), this is equivalent to that
% 1) in the exact sparse case ($q=0$), the sparsity $s$ is greater than
%a constant faction of the ambient dimension $p$;
% 2) under weak sparsity ($q \in(0,2)$),
the weak-$\ell_q$ radius satisfies
% \begin{equation}
% % s \gtrsim p \pth{\frac{r+c}{nh(\lambda)}}^{q/2}
$s \gtrsim p (\frac{r}{nh(\lambda)})^{q/2}$.
% \label{eqregPCA}
% \end{equation}
% for some constant $c \geq1$.
In the exact sparse case ($q=0$), %\eqref{eqregPCA}
this condition reduces to that the sparsity $s\asymp p$.
% $s$ is greater than a constant faction of the ambient dimension $p$.
% \nr{This remark seems to be repeating the preamble of Theorem 5. Need
%consolidation.}
% \label{rmkregPCA}
% \end{remark}

In the special case of $r=1$, a similar combinatorial procedure to
(\ref{eqBstar})--(\ref{eqcombi}) has been proposed in \cite{Vu12}.
Using Mendelson's results on empirical processes \cite{Mendelson10},
this procedure requires no sample splitting but can only be shown to
attain a convergence rate that is suboptimal in $\lambda$ \cite{Vu12},
Theorem 2.2: with $\lambda\to\infty$ and all the other parameters
fixed, the upper bound in \cite{Vu12} does not vanish. In contrast, the
optimal rate $\Psi$ decays at the rate $\lambda^{-(1-q/2)}$ when $k_q^*
< p$ and $\lambda^{-1}$ when $k_q^*=p$.
% \nr{Point out the suboptimality is when $\lambda$ is large?}
Comparing with the analysis in \cite{Vu12}, the proof of
Theorem~\ref{thmagg} is more elementary. By exploring the structure of
the difference between the sample covariance matrix and the true
covariance matrix, we obtained an upper bound that is optimal in all
parameters.
%here?} In fact, it appears that directly using results on empirical
%processes (such as \cite{Mendelson10} or Talagrand's theorem on
%second-order Gaussian chaos \cite{Talagrandchaining}, Theorem 2.5.2)
%yields a suboptimal dependence on both $\lambda$ and $r$.

%s3 #&#
\section{Adaptive estimation}\label{secadaptive}

The aggregation estimator constructed in Section~\ref{secUB} has been
shown to be rate optimal. However, it depends on the unknown parameters
and is computationally infeasible when $p$\vadjust{\goodbreak} is large. We construct in
this section an adaptive estimation procedure for principal subspaces
which is fully data driven and easily computable. Furthermore, it is
shown that the estimator attains the optimal rate of convergence
simultaneously over a large collection of the parameter spaces defined
in (\ref{paraspace}).

A key idea in our construction is a reduction scheme which reduces the
sparse PCA problem to a high-dimensional multivariate regression
problem. This method is potentially applicable to other sparsity
patterns of the leading eigenvectors. We first introduce the general
reduction scheme in Section~\ref{secreductionscheme} which transforms
the principal subspace estimation problem to a high-dimensional
multivariate regression problem. The specialization of this general
method under weak-$\ell_q$ constraint
%primary interest in this paper
will be detailed in Section~\ref{secreductionscheme2}.

%s3.1 #&#
\subsection{A general reduction scheme}\label{secreductionscheme}

The general reduction scheme involves four steps, which are introduced
in order below. The procedure used in step~2 for initial estimation
will be specified in Section~\ref{secreductionscheme2} for
weak-$\ell_q$ constrained parameter spaces. For ease of exposition, we
regard the rank $r$ as given in the statement below. Data-based choice
of $r$ will be discussed at the end of
Section~\ref{secreductionscheme2}.
\begin{longlist}
\item[\textit{Step} 1: \textit{Sample generation}.] Given the data
    matrix $\mathbf{X}$ in (\ref{eqmodel}) with $\sigma=1$, we generate
    an $n\times p$ random matrix $\tilde\mathbf{Z}$ with i.i.d.
    $N(0,1)$ entries which are independent of $\mathbf{U}$ and
    $\mathbf{Z}$, and form two samples $\mathbf{X}^i = \mathbf{X}+
    (-1)^i \tilde\mathbf{Z}$, $i=0,1$. Let $\mathbf{Z}^{i} = \mathbf
    {Z}+ (-1)^i \tilde\mathbf{Z}$ for $i=0,1$, then $\mathbf{Z}^0$ and
    $\mathbf{Z}^1$ are independent, and their entries are i.i.d.
    $N(0,2)$ distributed. Then the two samples $\mathbf{X}^0$ and
    $\mathbf{X}^1$ can be equivalently written as
%
%
%e30 #&#
\begin{equation}
\mathbf{X}^i = \mathbf{U}\mathbf{D}\mathbf{V}' +
\mathbf {Z}^i,\qquad i = 0,1. % \label{eqsplit}
\label{eqsamples}
\end{equation}
Let $\mathbf{S}^i = \frac{1}{n}(\mathbf{X}^i)'\mathbf{X}^i$,
$i=0,1$, be the sample covariance
matrices for the two samples.

\item[\textit{Step} 2: \textit{Initial estimation}.] We\vspace*{1pt}
    use the sample
    $\mathbf{X}^0$ to compute an initial estimator~$\mathbf{V}^0$. A
    specific procedure for computing the initial estimator
    $\mathbf{V}^0$ will be given in Section~\ref{secreductionscheme2}.
% This initial estimator $\VZERO$ can be obtained by ... \nr{We need to
%mention at least one specific method here.}

\item[\textit{Step} 3: \textit{Reduction to regression}.] Form
%
%
%e31 #&#
\begin{equation}
\label{eqregression-1} \bigl(\mathbf{X}^1\bigr)'
\mathbf{X}^0\mathbf{V}^0= \mathbf{V}\mathbf{A}+ \bigl(
\mathbf{Z}^1\bigr)'\mathbf{B},
\end{equation}
where $\mathbf{B}= \mathbf{X}^0\mathbf{V}^0$
%$ = \U\D\V'\VZERO+\ZZERO\VZERO$,
and
$\mathbf{A}= \mathbf{D}\mathbf{U}'\mathbf{B}$.
% $\bfA= \D\U'\U\D\V'\VZERO+\D\U'\ZZERO\VZERO$
% Note that conditioning on $\U$ and $\ZZERO$, both $\V\AA$ and $\BB$
%are fixed.
% Hence, \eqref{eqregression-1} becomes a regression problem,
%where $(\ZONE)'\bfB$ is an additive noise matrix with normal entries.
% However, since $\bfB$ does not have orthonormal columns, the noisy
%entries are not iid.
We now ``whiten'' the matrix in (\ref{eqregression-1}) as follows.
% , we introduce a further ``whitening'' step as follows.
Note that $\mathbf{B}= \mathbf{X}^0\mathbf{V}^0$ can be explicitly
computed after step~2.
Let its singular value decomposition be $\mathbf{B}= \mathbf
{L}\mathbf{C}\mathbf{R}'$, where
$\mathbf{L}\in\mathbb{R}^{n\times r}, \mathbf{C}\in\mathbb
{R}^{r\times r}$ and $\mathbf{R}\in
\mathbb{R}^{p\times r}$. Post-multiply both sides of (\ref{eqregression-1})
by $\frac{1}{\sqrt{2}}\mathbf{R}\mathbf{C}^{-1}$ to obtain
%
%
%e32 #&#
\begin{equation}
\label{eqregression} \mathbf{Y}= \bolds{\Theta}+ \mathbf{E},
\end{equation}
where
% \begin{equation}
% % \label{eqreg-terms}
$\mathbf{Y}= \frac{1}{\sqrt{2}} (\mathbf{X}^1)' \mathbf
{X}^0\mathbf{V}^0\mathbf{R}\mathbf{C}^{-1}$,
$\bolds{\Theta}= \frac{1}{\sqrt{2}} \mathbf{V}\mathbf{A}\mathbf
{R}\mathbf{C}^{-1}$ and
$\mathbf{E}= \frac{1}{\sqrt{2}}(\mathbf{Z}^1)'\mathbf{L}$.
% \end{equation}
% Compute the singular value decomposition of $\BB$ as $\BB= \bfL\bfD
% Let %$\bfC= \bfD\bfR$,
% $\bar{\AA} = \AA\bfR\bfD^{-1}$, $\bar{\BB} = \BB\bfR\bfD^{-1} =
% Form
% In the literature of high-dimensional regression,
% \eqref{eqrow-space} is usually referred to as the group sparsity
%constraint on the regression coefficients $\TT$.
We shall treat (\ref{eqregression}) as a regression problem, where
$\mathbf{Y}$ is the observed matrix, $\bolds{\Theta}$ is the signal
matrix of interest and $\mathbf{E}$ is the additive noise
matrix.\vadjust{\goodbreak}
Equivalently, we think of $\bolds{\Theta}$ as the coefficient matrix,
and the design matrix is $\mathbf{X}= \mathbf {I}_p$. The reason why
this is plausible will be detailed in
Section~\ref{secreductionscheme2}.

Given $\mathbf{Y}$, we propose the following method for computing
$\whT$. Define
% \nb{In this case, we cannot get $r\wedge(k-r)$!--ZM}
%
%
%e33 #&#
\begin{equation}
\label{eqt-k} t_k = r + \sqrt{2r \beta\log\frac{\mathrm{e}p}{k}} +
\beta\log\frac
{\mathrm{e}p}{k}.
\end{equation}
Fix an arbitrary $\delta\in(0,1)$. With slight abuse of notation, define
%
%
%e34 #&#
\begin{equation}
\label{eqpen} \operatorname{pen}(\bolds{\Theta}) = \operatorname {pen} \bigl(\bigl|
\operatorname{supp}(\bolds{\Theta})\bigr| \bigr),\qquad\mbox{where }
\operatorname{pen}(k) = (1+\delta)^2 \sum
_{i=1}^k t_i.
\end{equation}
Then the estimator for $\bolds{\Theta}$ is defined as
%
%
%e35 #&#
\begin{equation}
\label{eqtheta-est} \whT= \mathop{\operatorname{argmin}}_{\bolds{\Theta}\in\mathbb{R}^{p \times r}} \|\mathbf{Y}- \bolds{
\Theta}\|_{\mathrm{F}}^2 + \operatorname{pen}(\bolds{\Theta}).
\end{equation}
Such a penalized least squares approach has been widely used in
orthogonal regression with various choices of the penalty functions.
See, for example, Birg\'{e} and Massart \cite{Birge01} and Abramovich
et al. \cite{ABDJ06}.
\end{longlist}

%
%
%re4 #&#
\begin{remark}\label{rmkregression}
The penalized least squares estimator $\whT$ in (\ref{eqtheta-est}) can
be easily computed. Recall (\ref{eqregression}) and write the $i$th row
of the matrix $\mathbf {Y}$ as $\mathbf{y}_i$ and so $\mathbf{Y}=
[\mathbf{y}_1,\ldots,\mathbf {y}_p]'$. Let $\mathbf{y}_{(i)}$ denote
the row in $\mathbf{Y}$ with the $i$th largest $\ell_2$ norm, that is,
$\|\mathbf{y}_{(1)}\|\geq\|\mathbf{y}_{(2)}\|\geq\cdots\geq\|
\mathbf{y}_{(p)}\|$. Define
%
%
%e36 #&#
\begin{equation}
\label{eqkhat} \hat k =\mathop{\operatorname{argmin}}_{k \in[p]} \Biggl\{(1+
\delta)^2 \sum_{i=1}^k
t_i + \sum_{i=k+1}^p \|
\mathbf{y}_{(i)}\|^2 \Biggr\}.
\end{equation}
%
%Note that $\hat{k}$ is uniquely defined with probability $1$.
In case of multiple minimizers, $\hat{k}$ is chosen to be the smallest
one. It is also clear that $\hat k$ is easy to compute. With $\hat{k}$,
the estimator $\whT$ is given by $\widehat{\bolds{\Theta}}
=[\hat{\bolds{\theta}}_1,\ldots, \hat{\bolds{\theta}}_p]'$ where
\[
\hat{\bolds{\theta}}_i = \mathbf{y}_i \cdot{
\mathbf{1}_{ \{
{\|\mathbf{y}_i \|^2 > (1+\delta)^2 t_{\hat k}} \}}}.
\]
Note that $\hat{k}$ can be equivalently defined as $\operatorname
{argmin}_{k\in[p]} \sum_{i=1}^k [(1+\delta)^2 t_i -
\|\mathbf{y}_{(i)}\|^2]$. Therefore\vspace*{1pt}
$\|\mathbf{y}_{(\hat{k})}\|^2 > (1+\delta)^2 t_{\hat{k}}$ and
$\|\mathbf{y}_{(\hat{k}+1)}\|^2 \leq(1+\delta)^2 t_{\hat{k}+1}$. Since
$t_k$~is strictly\vspace*{-1pt} decreasing in $k$, we\vspace*{-2pt} obtain that
$\|\mathbf{y}_{(1)}\|^2 \geq \cdots\geq\|\mathbf{y}_{(\hat{k})}\|^2 >
(1+\delta)^2 t_{\hat{k}} \geq \|\mathbf{y}_{(\hat{k}+1)}\|^2
\geq\cdots.$ Thus, $|\operatorname {supp}(\whT)|=\hat{k}$.
\end{remark}

\textit{Step} 4: \textit{Final estimation}.
% In the final step, we find an estimator $\wh\TT$ for $\TT$ under model
Last but not least, we obtain the estimator $\whV$ for $\mathbf{V}$ by
orthonormalizing the columns of $\whT$. The orthonormalization can be
completed by the Gram--Schmidt procedure or QR factorization. The
estimated subspace is $\operatorname{span}(\whV) =
\operatorname{span}(\whT)$.

An important feature of the above reduction scheme is that the two
samples $\mathbf{X}^0$~and~$\mathbf{X}^1$ share the \emph{same}
realization of random factors $\mathbf{U}$ and their only difference is
in the noise matrices $\mathbf{Z}^0$ and $\mathbf{Z}^1$. This
is\vadjust{\goodbreak}
critical for maintaining the right level of signal-to-noise ratio in
the regression problem (\ref{eqregression}).
% when conditioning on $\U$ and $\ZZERO$.
In contrast, splitting the original sample into two halves as in
Section~\ref{secUB} does not achieve this goal here. Since our analysis
relies on the independence of $\mathbf{Z}^0$ and $\mathbf{Z}^1$, the
normality of
the noise is crucial to this scheme.

% \nb{---below is old---}
%
% Hereunder is the proposed procedure:
% \begin{enumerate}
% \item Use $\X_{(1)}$ to compute an initial estimator $\V_{(1)}$.
% % [with certain correct exclusion property]. So, w.h.p., $|\supp(
% \item
% \item Compute $\CC= \mathrm{diag}(\|\col{\BB}{1}\|,\dots, \|\col{
%Form
% \[
% \Y= \V\bar{\AA} + \Z_{(2)}\bar{\BB}.
% \]
% \item Compute an estimator $\wh{\TT}$ of $\TT= \V\bar{\AA}$ by the
%procedure \eqref{eqtheta-est} given later.
% \item Form the QR factorization of $\wh{\TT}$ as $\wh{\TT} = \wh\V\wh
% \end{enumerate}

%s3.2 #&#
\subsection{Sparse PCA and regression with group sparsity}\label{secreductionscheme2}

We now apply the general reduction scheme to the principal subspace
estimation problem with the parameter spaces defined in
(\ref{paraspace}). In what follows, we first introduce and study
a~\mbox{specific} estimator for the initial estimation step. Then we
derive properties of the proposed estimator for the regression with
group sparsity problem. Furthermore, we show that the general reduction
scheme paired with the two specific estimators leads to a final
estimator which adaptively achieves the optimal rates of estimation
over a large collection of the parameter spaces of interest. For
clarity of exposition, we regard the rank $r$ as given when introducing
the estimators. Data-driven choice of $r$ is discussed at the end of
this subsection.

\subsubsection*{Initial estimation}
Let $p_n\triangleq p\vee n$.
We construct the initial estimator $\mathbf{V}^0$ via the diagonal
thresholding method \cite{JohnstoneLu09} as follows:
% For the ease of exposition, we regard the rank $r$ as given.
% We discuss how to estimate it consistently from data is included at
%the end of this section.

%Recall that $\SZERO= \frac{1}{n}(\XZERO)'\XZERO$. $\VZERO$ is
%computed by the following steps:
%
\begin{longlist}[(3)]
\item[(1)] Define the set of features
%
%
%e37 #&#
\begin{equation}
\label{eqset-J} J = \bigl\{j\dvtx s^0_{jj} \geq2(1+ \alpha
\sqrt{\log{p_n}/n}) \bigr\},
\end{equation}
where $\{s^0_{jj}\}_{j=1}^p$ are the diagonal elements of $\mathbf{S}^0=
\frac{1}{n}(\mathbf{X}^0)'\mathbf{X}^0$, and $\alpha> 0$ is a
tuning parameter.

\item[(2)] Compute the first $r$ eigenvectors
$\{\hat{\mathbf{v}}^J_1,\ldots, \hat{\mathbf{v}}^J_r\}$ of
the submatrix $\mathbf{S}^0_{JJ}$.
% , denoted by $\wh\bfv^J_1,\dots, \wh\bfv^J_r$.

\item[(3)] Define $\mathbf{V}^0\in O(p,r)$, where
%
%
%e38 #&#
\begin{equation}
\label{eqinitial-V} \mathbf{V}^0_{J *} = \bigl[\hat{
\mathbf{v}}^J_1, \ldots, \hat{\mathbf{v}}^J_r
\bigr], \qquad\mathbf{V}^0_{J^c *} = \bolds{0}.
\end{equation}
\end{longlist}

The following result, proved in Section~7.5 in the supplementary
material \cite{appsm}, gives sufficient conditions on the model
parameters and the choice of $\alpha$ to guarantee that the initial
estimator $\mathbf{V}^0$ is reasonably close to $\mathbf{V}$, which
suffices for the initialization of our scheme.

%
%pr1 #&#
\begin{proposition}\label{propinitial-V}
Suppose that $\log n\geq M_0 \log\lambda$ for
some constant \mbox{$M_0 > 0$}. Suppose that
%% \log{\lambda}\leq M_0\log{p}, \quad\log{n}\leq M_0\log{p}
%
%
%e39 #&#
\begin{equation}
\label{eqcondition-initial1} n \bigl(\lambda^2 \wedge1 \bigr) \geq C_0
{(r+ \log{p})^2}/{\log{p}}
\end{equation}
and
%
%
%e40 #&#
\begin{equation}
\label{eqcondition-initial2} %\kappa^2 s \qth{\frac{\log(p\vee n)}{n(
s^2 \biggl(
\frac{\log(p \vee n)}{n\lambda^2} \biggr)^{1-q/2} < \kappa^{-2}
(2-q)^{q} / C_0
\end{equation}
for a sufficiently large constant $C_0>0$.
%and $n\geq C_0 \log p$ \nr{used to be $n\geq C_0 \log(p\vee n)$}.
If $\mathbf{V}^0$ is defined in (\ref{eqinitial-V}) with a sufficiently
large $\alpha\geq\sqrt{10(1+1/M_0)}$ in (\ref{eqset-J}),\vadjust{\goodbreak} then uniformly
over $\Theta= \Theta_q(s,p,\break r,\lambda)$, we have
%
%
%e41 #&#
\begin{equation}
\label{eqinitial-requirement} \bigl|\operatorname{supp}\bigl(\mathbf {V}^0\bigr)\bigr|\leq
k_q^*\quad\mbox{and}\quad \sigma_r \bigl(
\mathbf{V}'\mathbf{V}^0 \bigr)\geq1/2
\end{equation}
hold with probability at least $1 - C/[nh(\lambda)]$, where $k_q^*$ is
defined in (\ref{eqkqs}).
% , where $C$ is an absolute constant.
\end{proposition}

We note that condition (\ref{eqcondition-initial2}) is critical in
establishing the second claim in (\ref{eqinitial-requirement}), which
ensures that $\mathbf{V}^0$ is a reasonable estimator of $\mathbf {V}$.
Such a condition is needed for diagonal thresholding to work even when
$r = 1$. See, for example, condition~C3 in \cite{Paul05}, page~95.
Theorem 4.1 of \cite{Birnbaum12} showed that diagonal thresholding
could be suboptimal even under a stronger condition than
(\ref{eqcondition-initial2}).

% \nr{Add remark to condition \eqref{eqcondition-initial2} and
%compare to Paul's C3 in thesis. Analogy to the usual price to pay ($k$
%versus $k^2$) to use diagonal thresholding.}

%
%
%re5 #&#
\begin{remark}\label{rmkM0}
% \nb{Fix after revising conditions in prop 1! -- ZM}
% In Proposition~\ref{propinitial-V}, one possible choice of $\alpha$
%is $
When $M_0$ in Proposition~\ref{propinitial-V} is unknown, we replace it
by
%
%
%e42 #&#
\begin{equation}
\label{eqwhM0} \widehat{M}_0 = \log{n}/\log \bigl(
\sigma_1\bigl(\mathbf{S}^0\bigr)-2 \bigr),
\end{equation}
where $\sigma_1(\mathbf{S}^0)$ is the largest eigenvalue of $\mathbf
{S}^0$. This estimate works because \mbox{$\sigma_1(\mathbf{S}^0)-2$}
is an over-estimate of $\lambda$ with high probability
\cite{Paul07,Nadler08}, since the noise variance here is two. The
estimator (\ref{eqwhM0}) allows us to choose $\alpha$ in
(\ref{eqset-J}) without explicit knowledge of $M_0$.
\end{remark}

\subsubsection*{Orthogonal regression with group sparsity}
We first explain why we can treat~(\ref{eqregression}) as a regression
problem.
% Note that $\ZONE$ has iid $N(0,2)$ entries, and that the left
%singular vector matrix $\bfL$ of $\bfB$ has orthonormal columns.
% Conditioning on $\U$ and $\ZZERO$,
% the $\TT$ matrix in \eqref{eqregression} is fixed.
When we condition on the values of $\mathbf{U}$ and $\mathbf{Z}^0$, the
matrix $\mathbf{X}^0$ becomes deterministic. Thus, as deterministic
functions of $\mathbf{X}^0$, the matrices $\mathbf{V}^0, \mathbf{B},
\mathbf{L}, \mathbf{C}$ and $\mathbf{R}$ are also deterministic.
Furthermore, $\mathbf{A}$ and hence $\bolds{\Theta}$, as deterministic
functions of $\mathbf{U}$ and $\mathbf{B}$, are also deterministic. On
the other hand, $\mathbf{Z}^1$~is independent of both $\mathbf{U}$ and
$\mathbf {Z}^0$ and hence is independent of $\mathbf{X}^0$,
$\mathbf{B}$ and $\mathbf{L}$. Thus, the conditional distribution of
$\mathbf{Z}^1$ on $(\mathbf{U},\mathbf{Z}^0)$ always has i.i.d.
$N(0,2)$ entries, and so the conditional distribution of $\mathbf{E}$
has i.i.d. standard normal entries. Therefore, when we condition on the
values of $\mathbf{U}$ and $\mathbf{Z}^0$, problem (\ref{eqregression})
indeed reduces to a standard multivariate regression problem with
\emph{orthogonal design} and \emph{white noise}.

When the sparsity of $\mathbf{V}$ is specified as in (\ref
{paraspace}), we need
to consider the following parameter space for $\bolds{\Theta}$:
%
%
%e43 #&#
\begin{equation}
\label{eqrow-space} \mathcal{F}_q \bigl(s',p \bigr) = \bigl
\{ \bolds {\Theta}\dvtx\| \bolds{\Theta}\|_{q,w} \leq s'
\bigr\}
\end{equation}
%
% In this case, in the related regression problem \eqref{eqregression}
%denote by $\|\y\|_{[i]}$ and $\|\bolds{\theta}\|_{[i]}$ the $i$th
%largest row
%norm of the $\Y$ matrix and $\TT$ matrix respectively.
% The relevant parameter space for $\TT$ is $\calF_q(s',p)$ as defined
%in \eqref{eqrow-space}
with $q\in[0,2)$. The parameter $s'$ is typically different from $s$ in
(\ref{paraspace}), as it also depends on the other model parameters as
well as the realization of $\mathbf{U}$ and $\mathbf{Z}^0$. However,
this will not
cause any difficulty in practice, because the estimator proposed in
(\ref{eqtheta-est}) and the associated theorem below remain valid for
all values of $s'>0$.
% Moreover, $s'$ can be controlled with high probability.
In the literature of high-dimensional regression, (\ref{eqrow-space})
is usually referred to as the group sparsity constraint on the
regression coefficients $\bolds{\Theta}$.

For the estimator $\whT$ in (\ref{eqtheta-est}), we have following
upper bound on its risk. By the lower bounds in \cite{Lounici11} for $q
= 0$, the rates in Theorem \ref{thmseq} are optimal.

%
%
%th6 #&#
\begin{theorem}\label{thmseq}
Consider the regression problem %in \eqref{eqregression}
\[
\mathbf{Y}= \bolds{\Theta}+ \mathbf{E},
\]
where $\bolds{\Theta}\in\mathbb{R}^{p \times r}$ is deterministic
%is the $p\times r$ regression coefficients of interest
and $\mathbf{E}$ has i.i.d. $N(0,1)$ entries. Let the parameter space
$\mathcal{F}_q(s',p)$ be defined in (\ref{eqrow-space}) for some
$q\in[0,2)$ and $s'>0$. If $\beta> 2$ in~(\ref{eqt-k}) and
$\delta\in(0,1)$ in~(\ref{eqpen}), then there is a positive constant
$C$ that depends only on $q$, $\beta$ and $\delta$, such that the
estimator in (\ref{eqtheta-est}) satisfies
\[
\sup_{\bolds{\Theta}\in\mathcal{F}_q(s',p)}\mathsf{E}\|\whT- \bolds{\Theta}
\|_{\mathrm{F}}^2 \leq C k' \biggl(r + \log
\frac{\mathrm{e}p}{k'} \biggr),
\]
where
%
%
%e44 #&#
\begin{equation}
k' \triangleq\min \bigl\{k \in[p]\dvtx t_k^{q/2}
k \geq s' \bigr\} \label{eqkqs-z}
\end{equation}
for $t_k$ defined in (\ref{eqt-k}), and if the set in (\ref{eqkqs-z})
is empty, we set $k'=p$.
\end{theorem}

% \nr{Need to conform the definition of $k'$ and $k^*_q$. Use $\max$.
%
% $k \mapsto t_k^{q/2} k$ is increasing iff
% \[
% \frac{\beta\left(\frac{r}{\sqrt{2 \beta r \log\left(\frac{e p}{k}
% \]}

% For Theorem \ref{thmseq} to hold, there is no requirement on $(q,s')$
%as long as $q\in[0,2)$ and $s'>0$.

\subsubsection*{Adaptation}
With the above preparation, we are now ready to show that if we start
with a proper initial estimator $\mathbf{V}^0$ [such as that in
(\ref{eqinitial-V})] and estimate $\bolds{\Theta}$ by (\ref
{eqtheta-est}), then the estimator $\whV$ resulting from
orthonormalizing the columns of $\whT$ achieves the optimal rates of
convergence. We state the theorem in a slightly more general format. In
particular, it holds for the initial estimator in (\ref{eqinitial-V})
under the conditions of Proposition~\ref{propinitial-V}.
% \eqref{eqcondition-initial} and \eqref{eqcondition-initial2}.

%
%
%th7 #&#
\begin{theorem}[(Adaptation)]\label{thmupper-bound}
%Let $q \in[0,2)$.
Let $\lambda\geq C_0$ for some sufficiently large constant~$C_0$. Let
$\Theta= \Theta_q(s,p,r,\lambda)$ satisfy the conditions in
Theorem~\ref{thmagg}.
%sufficiently large absolute constant $M_0$},
Suppose that there exists an initial estimator $\mathbf{V}^0$
which\vspace*{1pt} satisfies (\ref{eqinitial-requirement}) with
probability at least $1 - C'/(nh(\lambda))$. Then the estimator $\whV$
obtained by orthonormalizing $\whT$ in (\ref{eqtheta-est}) with $\beta>
2$ in (\ref{eqt-k}) and $\delta\in(0,1)$ in (\ref{eqpen}) satisfies
% Suppose there exists an estimator $\VZERO$ such that with probability
%at least $1 - [nh(\lambda_1)]^{-1}$, $\supp|\VZERO|\leq k$ and $
%absolute constant. Then
%
\[
\sup_{\bolds{\Sigma}\in\Theta} \mathsf{E}\bigl\|\whV\whV' - \mathbf{V}
\mathbf{V}'\bigr\|_{\mathrm{F}}^2 \leq 2 \bigl(r
\wedge(p-r) \bigr) \wedge C\Psi \bigl(k^*_q,p,r, n,\lambda \bigr),
%C \left[r \wedge\frac{ k_q^*}{nh(
%p}{k_q^*}\right)\right].
\]
where $k_q^*$ is defined in (\ref{eqkqs}), and $C>0$ is a constant
depending only on $q,\beta$ and $\delta$.
%$q,\kappa,\beta$ and $\delta$.
\end{theorem}

%We note that \eqref{eqassumption-superbigsnr} is needed in
%establishing the general upper bounds via aggregation in
%Theorem~\ref{thmagg}.
%In addition,
We note that the assumption $\lambda> C_0$ is imposed to ensure that
% $\BB$ is of rank $r$, and hence
the ``whitening'' procedure in step~3 of the reduction scheme can be
performed.
%aggregation upper bounds??}

It is interesting to compare the statement of
Theorem~\ref{thmupper-bound} to the minimax lower bound in Theorems
\ref{thmlower-bound}--\ref{thmlower-bound-0} as well as the
performance of the combinatorial aggregation estimator
$\widehat{\mathbf{V}}^*$ established in Theorem~\ref{thmagg}. For any
parameter space $\Theta= \Theta_q(s,p,r,\lambda)$ such that the\vspace*{1pt}
conditions of Proposition~\ref{propinitial-V}
% the conditions \eqref{eqcondition-initial} and
hold, we could use the $\mathbf{V}^0$ in (\ref{eqinitial-V}), and the
resulting $\whV$ is guaranteed to achieve the optimal rates of
convergence on $\Theta$, which matches the performance of the
aggregation estimator for any $q>0$. Moreover, in this case both
$\mathbf{V}^0$ and $\whV$ can be efficiently computed. Hence $\whV$ can
be used in practice while $\whV^*$ is computationally
intensive.\vadjust{\goodbreak}
However, in the exact sparse case of $q=0$, the upper bound in
Theorem~\ref{thmupper-bound} depends on the rank $r$ linearly through
$s r$, while the true minimax rate in Theorem~\ref{thmlower-bound-0}
depends on $r$ quadratically through $r(s-r)$, which is smaller than
$rs$ if $s-r$ is small. The suboptimality of $\whV$ in this specific
regime is partially due to the fact that our reduction scheme
transforms the problem into a regression problem without taking account
of the orthogonality structure of the parameter space.
% By \cite{Birnbaum12} and \cite{Ma11}, a sufficient condition for the
%condition of Theorem \ref{thmupper-bound} to hold is
% \begin{equation}
% \frac{k^2\log{p}}{n\lambda_r^2} \leq\eps
% \end{equation}
% for a sufficiently small constant $\eps> 0$.
% \nb{Remark on the condition of Theorem \ref{thmupper-bound} to be
%added! One could give some sufficient condition on it!}

%
%
%re6 #&#
\begin{remark}
Theorem~\ref{thmupper-bound} shows that any estimator $\mathbf{V}^0$
satisfying (\ref{eqinitial-requirement}) can be used to produce an
adaptive estimator. So the task of constructing adaptive optimal
estimators is reduced to constructing a ``reasonable'' estimator.
% In fact, the requirement of $|\supp(\VZERO)|\leq k^*$ in
%to generate one more sample $\XTWO$ for constructing the initial
%estimator,
% \nb{start here!!!!}
\end{remark}

% \begin{remark}
% \nb{Give more remarks on the construction of the two samples here!
%Note that this is not equal to split sample, because we need to have
%the same realization of the random effects $\U$ in $\XONE$ and $
%to an optimal estimator.}
% \end{remark}

\subsubsection*{Consistent estimator of $r$}
Last but not least, we discuss how to construct a~consistent estimator
of $r$ based on data. To this end, recall the definition of the set $J$
in (\ref{eqset-J}), and the matrix $\mathbf{S}^0_{JJ}$. We propose to
estimate $r$ by
%
%
%e45 #&#
\begin{equation}
\label{eqr-hat} \hat{r} = \max \bigl\{l\dvtx\sigma_l\bigl(
\mathbf{S}^0_{JJ}\bigr) > 2(1+\delta_{|J|}) \bigr
\},
\end{equation}
where for any $m>0$ and $M_0$ in the conditions of
Proposition~\ref{propinitial-V},
% \eqref{eqcondition-initial}
we define
\[
\delta_m = 2(\sqrt{m / n}+t_m) + (\sqrt{m/
n}+t_m)^2
\]
with
$t_m^2 =
\frac{2}{n}((m+1)\log(\mathrm{e}p)+ (1+2/M_0)\log{n})$.
% $t_m^2 = \frac{2(1+M_0)\log{p}}{n} + \frac{2\log{n}}{n} + \frac{2m(
Here, we regard $M_0$ as known. Otherwise, we could always replace it
with the estimator (\ref{eqwhM0}) proposed in Remark~\ref{rmkM0}. Note
that the estimator (\ref{eqr-hat}) could be easily integrated with the
diagonal thresholding method for computing $\mathbf{V}^0$. In particular,
$\hat{r}$ can be computed after we select the set $J$
in~(\ref{eqset-J}).

For this estimator, we have the following result.
%
%
%pr2 #&#
\begin{proposition}
\label{propr-hat} Under the condition of
Proposition~\ref{propinitial-V}, $\hat{r} = r$ holds with probability
at least $1-C [nh(\lambda)]^{-1}$.
\end{proposition}
Under
% the conditions \eqref{eqcondition-initial} --
the conditions of Proposition~\ref{propinitial-V} and
Theorem~\ref{thmupper-bound}, Proposition~\ref{propr-hat} implies that
the conclusion in Theorem~\ref{thmupper-bound} still holds if we
replace $r$ by $\hat{r}$.

% \subsection{Orthogonal Regression with Group Sparsity}

% This subsection is devoted to analyzing the fourth step of the
%adaptive algorithm proposed in the previous subsection. To this end,
%we investigate here optimal rates (mainly the upper bound part) for
%estimation in orthogonal regression model with group sparsity.
%
% Consider the model
% \begin{equation}
% \label{eqgroup-model}
% \Y= \TT+ \Z,
% \end{equation}
% where $\Y= [\y_1,\dots,\y_p]'$, $\TT= [\bolds{\theta}_1,\dots,\bolds{
%= [\z_1,\dots, \z_p]'$ are all $p\times r$ matrices. $\Y$ is the
%observed responses, and $\Z$ has iid $N(0,1)$ entries.

%s4 #&#
\section{Numerical experiments}\label{secnumeric}
In this section, we report simulation results comparing the adaptive
method proposed in Section~\ref{secadaptive} with the iterative
thresholding method proposed in Ma \cite{Ma11}.

In all the results reported here, the sample size $n = 1000$ and the
ambient dimension $p = 2000$.
We focus on the case of exact sparsity, that is, $q = 0$.
The sparsity parameter $s$ takes value in $\{40, 80, 120, 160, 200\}$,
% $9$ equally-spaced integer values between $40$ and $200$.
% For each value of $s$,
and the rank $r$ takes value in $\{1, 5, 10, 20\}$. For each $(s,r)$
combination, the $\mathbf{V}$ matrix is obtained from orthonormalizing
an $p\times r$ matrix $\mathbf{M}$ where $\mathbf{M}_{i *}$ have i.i.d.
$N(0,i^4)$ entries for $i=1,\ldots, s$
% \nr{need to explain why: to fulfill (approximately?) the weak-lq
%constraint. what is the value of $q$ we are fitting here? -- yih}
and $\mathbf{M}_{i *}=0$ for all $i>s$.
We set the variances of different rows to be different so that the
ordered norms of the nonzero rows in $\mathbf{V}$ also exhibit fast decay.
When $r = 1$, the spike size $\lambda_1 = 20$.
When $r > 1$, the $\lambda_i$'s take $r$ equispaced values such that
$\lambda_r = 10$ and $\lambda_1 = 20$.

When implementing the method in Section~\ref{secadaptive}, we take
$\alpha= 3$ in~(\ref{eqset-J}), $\beta= 2.1$ in~(\ref{eqt-k}) and
$\delta= 0.05$ in~(\ref{eqpen}) in all the simulations reported here.
In addition, we made a slight modification to the proposed method to
obtain better numerical results. We first run the method to obtain an
estimator, denoted by $\whV_1$. Then we switch the roles of
$\mathbf{X}^0$ and $\mathbf{X}^1$ and run the proposed procedure again
to obtain a~second estimator~$\whV_2$. Finally, we use the $r$ leading
eigenvectors of $\whV_1\whV_1'+\widehat{\mathbf{V}}_2\whV_2'$ as the
columns of the final estimator $\whV$. By Theorem 10 in Section~7.11 in
the supplementary material~\cite{appsm}, we have
\[
\bigl\|\whV\whV' - \mathbf{V}\mathbf{V}\bigr\| _{\mathrm{F}} \leq\bigl\|
\whV_1\whV_1' + \whV_2
\whV_2' - 2\mathbf{V}\mathbf {V}'
\bigr\|_{\mathrm{F}} \leq \sum_{i=1,2} \bigl\|
\whV_i\whV_i' - \mathbf{V}\mathbf {V}
\bigr\|_{\mathrm{F}}.
\]
Here, the first inequality holds because $\sigma_r(\whV_1\whV_1' +
\whV_2\whV_2') \geq\sigma_r(\widehat{\mathbf{V}}_1\whV_1') = 1$ and
$\sigma_{r+1}(2\mathbf{V}\mathbf{V}') = 0$, while the second is by the
triangle inequality. By the last display, the theoretical results in
Section~\ref{secadaptive}, which apply to both $\whV_1$ and $\whV_2$,
also apply to the final estimator $\whV$.
% {\blue[Do we need to spell out the arguments here?]}
When implementing the iterative thresholding method in Ma~\cite{Ma11}, we
set all tuning parameters at their recommended values.

%
%
%t1 #&#
\begin{table}%[b]
\caption{Average loss $\|\whV\whV' -
{\mathbf{V}}{\mathbf{V}}'\|^2_{\mathrm{F}}$ over $50$ repetitions for
each $(s,r)$ combinations}\label{tabnumeric}
\tabcolsep=0pt
\begin{tabular*}{\tablewidth}{@{\extracolsep{\fill}}@{}lcccccc@{}}
\hline
&& \multicolumn{5}{c@{}}{$\bolds{s}$}\\[-6pt]
&& \multicolumn{5}{c@{}}{\hrulefill}\\
$\bolds{r}$ & \textbf{Method} & \textbf{40} & \textbf{80} & \textbf
{120} & \textbf{160} & \textbf{200}\\
\hline
\phantom{0}1 & RegSPCA & 0.0236 & 0.0660 & 0.0892 & 0.1074 & 0.1754 \\
& ITSPCA & 0.0117 & 0.0366 & 0.0483 & 0.0619 & 0.0712
\\[3pt]
\phantom{0}5 & RegSPCA & 0.0348 & 0.0718 & 0.1134 & 0.1470 & 0.1992 \\
& ITSPCA & 0.0520 & 0.1209 & 0.1848 & 0.2368 & 0.3042
\\[3pt]
10 & RegSPCA & 0.0544 & 0.1247 & 0.1777 & 0.2394 & 0.3052 \\
& ITSPCA & 0.0914 & 0.2284 & 0.3535 & 0.4866 & 0.6313
\\[3pt]
20 & RegSPCA & 0.0640 & 0.1826 & 0.2904 & 0.4030 & 0.5083 \\
& ITSPCA & 0.1185 & 0.3740 & 0.6449 & 0.9045 & 1.1715\\
\hline
\end{tabular*}
\end{table}

Table~\ref{tabnumeric} summarizes the average squared Frobenius losses
of the proposed method (RegSPCA) and the iterative thresholding method
(ITSPCA) over $50$ repetitions for each $(s,r)$ combination.
Table~\ref{tabnumeric} shows that for all values of the sparsity
parameter, RegSPCA outperformed ITSPCA when $r = 5, 10$ or $20$, while
ITSPCA led to smaller average losses when $r = 1$. This demonstrates
the competitiveness of RegSPCA in the group sparse setting considered
in the present paper. On the other hand, we note that ITSPCA was not
designed specifically for handling the group sparsity structure which
is the case when $r>1$, and hence its underperformance is not
unexpected.

%%%\iffalse
%%%\begin{tabular}{|c|cc|cc|cc|cc|}
%%%\hline
%%%$s$ & \multicolumn{2}{|c|}{$r=1$} & \multicolumn{2}{|c|}{$r=5$} &
%%%\multicolumn{2}{|c|}{$r=10$} & \multicolumn{2}{|c|}{$r=20$}\\
%%%\hline
%%%& RegSPCA & ITSPCA & RegSPCA & ITSPCA & RegSPCA & ITSPCA & RegSPCA &
%%%ITSPCA\\
%%%\hline
%%%$40$ & $0.0236$ & $0.0117$ & $0.0348$ & $0.0520$ & $0.0544$ &
%$0.0914$
%%%& $0.0640$ & $0.1185$ \\
%%%% $23.7600$ & $22.8800$ & $31.0800$ & $35.5800$ & $32.3800$ &
%$40.3000$
%%%%& $34.6800$ & $41.5000$ \\
%%%$80$ & $0.0660$ & $0.0366$ & $0.0718$ & $0.1209$ & $0.1247$ &
%$0.2284$
%%%& $0.1826$ & $0.3740$ \\
%%%% $34.6400$ & $37.6600$ & $56.6600$ & $61.3800$ & $59.5400$ &
%$70.1400$
%%%%& $64.9400$ & $80.7800$ \\
%%%$120$ & $0.0892$ & $0.0483$ & $0.1134$ & $0.1848$ & $0.1777$ &
%$0.3535$
%%%& $0.2904$ & $0.6449$ \\
%%%% $44.5000$ & $50.8200$ & $80.4600$ & $87.4600$ & $86.5400$ &
%$96.6000$
%%%%& $93.8600$ & $112.8600$ \\
%%%$160$ & $0.1074$ & $0.0619$ & $0.1470$ & $0.2368$ & $0.2394$ &
%$0.4866$
%%%& $0.4030$ & $0.9045$ \\
%%%% $50.8600$ & $60.2400$ & $105.9200$ & $113.2000$ & $113.3200$ &
%%%%$123.0600$ & $119.6200$ & $138.9600$ \\
%%%$200$ & $0.1754$ & $0.0712$ & $0.1992$ & $0.3042$ & $0.3052$ &
%$0.6313$
%%%& $0.5083$ & $1.1715$ \\
%%%% $59.6400$ & $81.9400$ & $120.6800$ & $130.5000$ & $137.1600$ &
%%%%$146.3200$ & $145.8000$ & $164.6600$
%%%\hline
%%%\end{tabular}
%%%
%%%\fi

%s5 #&#
\section{Discussions}\label{secdiscussion}

We have focused in the present paper on the estimation of the principal
subspace $\operatorname{span}(\mathbf{V})$ under the loss (\ref
{eqloss}). The minimax rates of
convergence are established and a computationally efficient adaptive
estimator is constructed.

Both the current paper and Ma \cite{Ma11} consider the problem of sparse
subspace estimation under the spiked model, but they differ in several
important ways. First, in addition to the sparsity constraint on the
leading eigenvectors, the current paper requires them to share support.
This extra assumption is motivated by real data applications. For
instance, if the observed vectors are the leading Fourier coefficients
of random functions with a common covariance kernel, then we expect the
leading eigenvectors to have large coefficients only at low frequency
coordinates so that the resulting leading eigen-functions in the time
domain are smooth. Second, Ma \cite{Ma11} focused on the error upper
bounds of an adaptive estimator with the subspace rank $r$ assumed to
be a fixed constant. Whether the dependence of the bounds on $r$ is
optimal was not studied. The current paper conducts an investigation on
the dependence of the minimax rates on key model parameters, including
$r$ which can grow with $n$ and $p$. Last but not least, we have
focused exclusively on the subspace $\operatorname{span}(\mathbf{V})$
which is natural when the
spikes are of the same order, while Ma \cite{Ma11} considered estimating
subspaces spanned by the first few rather than all columns of $\mathbf
{V}$. The
optimal rates of the latter estimation problem is of most interest when
the spikes scale at different rates with $n$ and $p$, which we leave as
an interesting problem for future research.

A problem closely related to principal subspace estimation is the
estimation of the whole covariance matrix $\bolds{\Sigma}$ under the same
structural assumption (\ref{paraspace}).
%In this case it is more natural to use the spectral norm as the loss
%function $L(\wh\bfSigma, \bfSigma) = \| \wh\bfSigma- \bfSigma\|^2$.
Both minimax estimation and adaptive estimation are of significant interest.
% Another relevant question is whether a plug-in estimator of the type $
%estimates of $ \LL$ and $\sigma^2$ respectively, can be rate optimal
%under the spectral norm loss.
Results on minimax rates under the spectral norm loss
$L(\widehat{\bolds{\Sigma}}, \bolds{\Sigma}) = \|
\widehat{\bolds{\Sigma}}- \bolds{\Sigma}\|^2$ can be found in~\cite{CMW13}.

It is interesting to extend the aggregation method in
Section~\ref{secUB} to other settings beyond sparsity or weak $\ell_q$
constraints. In the exact sparse case ($q=0$), note that the
rate-optimal estimator in (\ref{eqagg}) is constructed by choosing the
best estimator from a collection of estimators, each of which is
designed for a specific sparsity pattern. Theorem~\ref{thmagg} can now
be interpreted as an oracle inequality for the average risk, which is
within a constant factor\vspace*{-2pt} of the oracle risk $\frac{r(k-r)}{n
h(\lambda)}$ plus the excess risk $\frac{1}{n h(\lambda)} \log
{p\choose k}$.
%Apart from this model-selection idea, one could also consider
%aggregation by means of convex
One immediate generalization of Theorem~\ref{thmagg} is that we can
also construct aggregated estimators if it is known that the true
principle subspace belongs to a collection of $N$ subspaces. Then the
excess risk does not exceed $\frac{1}{n h(\lambda)} \log N$.
% It is much more challenging to obtain an oracle inequality with the
%constant one, which implies an upper bound on the \emph{minimax
%regret}. The current aggregation method based on the equal sample
%splitting is, however, not sufficient to achieve this goal.

%we want to state Theorem~\ref{thmagg} in terms of minimax regret by
%giving an oracle inequality with constant one, and maybe some
%high-probability upper bound. Moreover, we should formulate the
%statement in terms of \ul{union of $N$ subspaces} and assess the
%minimax regret involving $\log N$.
%Let $\ntok{V_1}{V_N}$ be subspaces of $\reals^p$ with dimension no
%more than $k$. Then
%Let $\hV_*$ be constructed similarly
%
%Generalization to union of subspaces:
% Let $\ntok{V_1}{V_N}$ be subspaces of $\reals^p$ with dimension no
%more than $k$. Then
%% such that $\max_{i\in[N]} \dim V_i = k$.
% \begin{equation}
% \label{eqaggsub}
% \label{thmaggsub}

It should be noted that our analysis in this paper relies on the
normality of the model, which allows us to express the sample in the
form of (\ref{eqmodel}). In particular the adaptive procedure requires
the independence of $\mathbf{Z}^0$ and $\mathbf{Z}^1$, which is a
consequence of the
normality of the noise. It is unclear whether the same results hold for
all noise distributions with sub-Gaussian tails. It is an interesting
problem to study the robustness of the adaptive procedure and to extend
the results to other noise distributions.

% The adaptation procedure proposed in the current paper shows that
%sparse PCA is connected to the Gaussian sequence model.
% Moreover, the optimal rates for the sparse PCA problem derived in the
%present paper coincide with those for the regression problem in
%theoretical question is whether certain forms of the two problems are
%indeed asymptotically equivalent to each other in the Le Cam's sense
% Such an asymptotic equivalence result would enable deeper
%understanding of the sparse PCA problem, and guide the development of
%other adaptive estimation procedures by borrowing the insights from
%the regression problem.

% \item Relationship to regression problem (general device)
% \item Particularization

%s6 #&#
\section{Proofs}\label{secproof}
In this section we prove Theorems
\ref{thmlower-bound-0}, \ref{thmagg} and \ref{thmupper-bound}. The
proofs of the other results, together with those of the key lemmas and
some additional technical arguments, are given in the supplementary material~\cite{appsm}.
% the main theorems and the key technical lemmas. Proofs of
%Propositions \ref{propinitial-V} and \ref{propr-hat}, together with
%some additional technical arguments, are given in the appendix.
% \nb{double check after restructuring!!!!}

%s6.1 #&#
\subsection{\texorpdfstring{Proof of Theorem \protect\ref{thmlower-bound-0}}
{Proof of Theorem 3}}\label
{secpflower-bound}

%%%\iffalse
%%%%
%%%\begin{remark}
%%%\nr{Temporary: discussion on assumptions}
%%%%
%%%\begin{itemize}
%%%%
%%%\item\eqref{eqassumption-lambdan}: needed in Wishart tail bound,
%%% Used in oracle upper bound as well as the aggregation upper
%%% bound on $\U$. Consequence of
%%% Proposition~\ref{propwishart-bd}.
%%%
%%%This is mild because all polynomial decay in $n$ is captured.
%%%Superpolynomial in $n$ is however impossible.
%%%% \item\eqref{eqassumption-rn}: need in Wishart tail bound. Used
%%%%in aggregaation upper bound on $\U$. Consequence of
%%%%Proposition~\ref{propwishart-bd}.
%%%
%%%\item\eqref{eqassumption-bigsnr}:
%%%
%%%\item\eqref{eqassumption-n-k-lambda}:
%%%
%%%\item\eqref{eqassumption-prk}:
%%%
%%%\item\eqref{eqassumption-superbigsnr}:
%%%
%%%\item\nb{$k^*\log(ep/k^*)/nh(\lambda)\leq c$ and $\log n\lambda
%%% /(n\lambda)\leq c$}: used in the proof of
%%% Theorem~\ref{thmupper-bound}.
%%%\end{itemize}
%%%%
%%%\end{remark}
%%%%
%%%\fi

% Two methods
% \begin{enumerate}
% \item Lower bound by more informative model (equivalent to
%regression): optimal except when $\lambda$ is small.
% \item Two lower bound:
% \begin{enumerate}
% \item by oracle risk if we know the support. This is clear.
% \item by the risk of estimating $v_1$, the singular vector for $1+
%loss?

%Starting from the exact sparse case ($q=0$),
We first give a lower bound on the oracle risk where we know beforehand
the row support of $\mathbf{V}$. This corresponds to a $k$-dimensional
unstructured PCA problem, where the goal is to estimate the $r$ leading
singular vectors of the covariance matrix. In view of the upper bound
in Theorem~\ref{thmoracle-ub}, the rates are minimax optimal.

%
%
%th8 #&#
\begin{theorem}[(Oracle risk: lower bound)]\label{thmoracle-lb}
% Let $p=k$ and $r \in[k]$.
Let $\Theta= \Theta_0(k,k,r,\lambda)$.
% where $\lambda> 0$.
Then
%
%
%e46 #&#
\begin{equation}
\inf_{\whV}\sup_{\bolds{\Sigma}\in\Theta} \mathsf {E}\bigl\|\whV
\whV' - \mathbf{V}\mathbf {V}'\bigr\|_{\mathrm{F}}^2
\gtrsim r \wedge(k-r) \wedge\frac{r(k-r)}{nh(\lambda)}. \label{eqoracle}
\end{equation}
%
% where $c$ is an absolute constant. % depending only on $\kappa$.
\end{theorem}

To prove Theorem~\ref{thmoracle-lb}, we use a minimax lower bound due
to Yang and Barron~\cite{YB99}, Section~7, via \emph{local} metric
entropy, which in turn relies on an argument by Birg{\'e}
\cite{Birge83}.
% The situation here is slightly different from that in \cite{YB99} in
%the sense that we use global covering number instead of packing number
%to derive bounds on local packing number.
For completeness, we state the result in Proposition~\ref{propyb} and
provide a~short proof in Section~7.8 in the supplementary
material~\cite{appsm}.
% The method of local metric entropy in an $\frac{1}{
The method of local metric entropy in an
$\frac{1}{\sqrt{n}}$-neighborhood dates back to Le Cam \cite{LeCam73}.
The advantage of this method is that it only relies on the analytical
behavior of the metric entropy of the parameter space, thus allowing us
to sidestep constructing explicit packing set in the parameter space.

%
%pr3 #&#
\begin{proposition}\label{propyb}
Let $(\Theta,\rho)$ be a totally bounded metric space and
$\{P_\theta\dvtx\break \theta\in\Theta\}$ a collection of probability
measures. For any $E \subset\Theta$, denote by $\mathcal
{N}(E,\epsilon)$ the
$\epsilon$-covering number of $E$, that is, the minimal number of balls
of radius $\epsilon$ whose union contains $E$. Denote by
$\mathcal{M}(E,\epsilon)$ the $\epsilon$-packing number of $E$, that
is, the
maximal number of points in $E$ whose pairwise distance is at least
$\epsilon$. Put
%
%
%e47 #&#
\begin{equation}
\label{eqKLd} A \triangleq\sup_{\theta\neq\theta'} \frac{D(P_\theta||
P_{\theta'})}{\rho^2(\theta,\theta')}.
\end{equation}
If there exist $0 < c_0 < c_1 < \infty$ and $d \geq1$ such that
%
%
%e48 #&#
\begin{equation}
\label{eqNepsilon} \biggl( \frac{c_0}{\epsilon} \biggr)^{d} \leq\mathcal{N}(
\Theta,\epsilon) \leq \biggl( \frac{c_1}{\epsilon} \biggr)^{d}
\end{equation}
for all $0 < \epsilon< \epsilon_0$. Then
%, where $\epsilon_0^2 = \frac{\log2}{24} \frac{d c_0^2}{A c_1^2}$.
%Then
%
%
%e49 #&#
\begin{equation}
\label{eqyb} \inf_{\hat{\theta}} \sup_{\theta\in\Theta}
\mathsf{E}_{\theta} \bigl[\rho^2 \bigl(\hat{\theta}(X),\theta
\bigr) \bigr] \geq \frac{c_0^2}{840 c_1^2} \biggl( \frac{d }{A} \wedge
\epsilon_0^2 \biggr).
\end{equation}
\end{proposition}
%
%We conclude this subsection by proving Proposition~\ref{propyb}:

We also need the following result regarding the metric entropy of the
Grassmannian manifold $G(k,r)$ due to Szarek \cite{Szarek82}.

%
%le1 #&#
\begin{lemma}\label{lmmszarek}
For any $\mathbf{V}\in O(k,r)$, identifying the subspace
$\operatorname{span}(\mathbf{V})$ with its projection matrix
$\mathbf{V}\mathbf{V}'$, define the metric on $G(k,r)$ by
$\rho(\mathbf{V}\mathbf{V} ',\mathbf{U}\mathbf{U}') =
\|\mathbf{V}\mathbf{V}'-\mathbf {U}\mathbf{U}'\|_{\mathrm{F}}$. Then
for any $\epsilon\in(0, \sqrt{2 (r\wedge(k-r))}]$,
%
%
%e50 #&#
\begin{equation}
\label{eqszarek} \biggl( \frac{c_0}{\epsilon} \biggr)^{r(k-r)} \leq\mathcal{N}
\bigl(G(k,r),\epsilon \bigr) \leq \biggl( \frac{c_1}{\epsilon} \biggr)^{r(k-r)},
\end{equation}
where $c_0,c_1$ are absolute constants. Moreover, for any $\mathbf
{V}\in
O(k,r)$ and any \mbox{$\alpha\in(0,1)$},
%
%
%e51 #&#
\begin{equation}
\mathcal{M} \bigl(B(\mathbf{V},\epsilon), \alpha\epsilon \bigr) \geq \biggl(
\frac{c_0}{\alpha c_1} \biggr)^{r(k-r)} \label{eqlocalGkr}.
\end{equation}
\end{lemma}

\begin{pf}
Note that $\rho(\mathbf{V}\mathbf{V}',\mathbf{U}\mathbf{U}') = \sqrt{2}
\|(\mathbf{I}- \mathbf{V}\mathbf{V}') \mathbf{U}\mathbf
{U}'\|_{\mathrm{F}}$, in view of (\ref{eqlossp}). This metric is
unitarily invariant; see $\rho_{\alpha}'$ in \cite{Szarek82}, Remark~5,
page~175. Applying \cite{Szarek82}, Proposition 8, page 169, with
$\alpha(\cdot) = \Vert{\cdot} \Vert$ gives (\ref{eqszarek}). By the
proof of equation (158) in the supplementary material~\cite{appsm}
% \eqref{eqMBeps},
for any $\epsilon\in(0, \sqrt{2 (r\wedge(k-r))}]$ and any $\alpha
\in(0,1)$, there exists $\mathbf{V}^* \in O(k,r)$ such that $\mathcal
{M}(B(\mathbf{V}
^*,\epsilon), \alpha\epsilon) \geq(\frac{c_0}{\alpha c_1})^{r(k-r)}$.
Now for any $\mathbf{V}\in O(k,r)$, there exists $\mathbf{T}\in
O(p)$, such that $\mathbf{V}=
\mathbf{T}\mathbf{V}^*$. Then (\ref{eqlocalGkr}) holds since the
metric $d$ is
unitarily invariant.
\end{pf}

\begin{pf*}{Proof of Theorem~\ref{thmoracle-lb}}
For the purpose of lower bound, we consider the special case of
$\lambda
_1 = \cdots=\lambda_r = \lambda$, that is, $\bolds{\Sigma}= \lambda
\mathbf{V}\mathbf{V}' +
\mathbf{I}_k$.
%A simple calculation of the Kullback-Leibler divergence yields
Note that the Kullback--Leibler divergence between normal distributions
is given by $D(N(0, \bolds{\Sigma}_1) || N(0, \bolds{\Sigma}_0)) =
\frac{1}{2}(\operatorname{Tr}(\bolds{\Sigma}_0^{-1} \bolds{\Sigma
}_1 -
\mathbf{I}_k) - \log\det\bolds{\Sigma}_0^{-1}\bolds{\Sigma}_1)$.
Then for any $\mathbf{U},\mathbf{V}\in O(k,r)$,
we have
% D(N(0, \lambda\V\V' + \I_k)^n || N(0, \lambda\U\U' + \I_k)^n) = n
%h(\lambda) \fnorm{\V\V'-\U\U'}^2.
% \label{eqKLA}
%
%
%e52 #&#
\begin{eqnarray}
\label{eqKLA} && D \bigl(N \bigl(0, \lambda\mathbf{V}\mathbf{V}' +
\mathbf{I}_k \bigr)^n || N \bigl(0, \lambda\mathbf{U}
\mathbf{U}' + \mathbf{I}_k \bigr)^n \bigr)
\nonumber
\\
&&\qquad= \frac{n}{2} \operatorname{Tr} \biggl( - \frac{\lambda
}{\lambda+1}
\mathbf{V}\mathbf{V}' + \lambda\mathbf{U}\mathbf{U}' -
\frac{\lambda^2}{\lambda+1} \mathbf{V}\mathbf{V}'\mathbf{U}
\mathbf{U}' \biggr)
\\
&&\qquad= \frac{n \lambda^2}{2(\lambda+1)} \bigl(r - \bigl\|\mathbf {U}'\mathbf{V}
\bigr\|_{\mathrm{F}}^2 \bigr)= \frac{n
h(\lambda)}{2} \bigl\|\mathbf{V}
\mathbf{V}'-\mathbf{U}\mathbf{U}'\bigr\|
_{\mathrm{F}}^2,\nonumber
\end{eqnarray}
where the first and second inequalities are by the matrix inversion
lemma and the fact that $\operatorname{Tr}(\mathbf{V}\mathbf
{V}')=\operatorname{Tr}(\mathbf{V}'\mathbf{V})=r$, respectively. In
view of (\ref{eqKLd}), we have $A = n h(\lambda)/2$. Applying
Proposition~\ref{propyb} with $\epsilon_0 = \sqrt{2 (r \wedge(k-r))}$
yields the desired~(\ref{eqoracle}).
\end{pf*}

\begin{pf*}{Proof of Theorem~\ref{thmlower-bound-0}}
Let $\Theta= \Theta_0(s,p,r,\lambda)$. By definition (\ref{eqkqs}),
$k^*_0$ coincides with $s$.
% Under the assumption of Theorem~\ref{thmlower-bound-0}, our goal is to
%prove the following nonasymptotic lower bound: if $r \leq p + 1 - s$,
%then
% %for any $1 \leq r \leq k \leq p$,
% \begin{equation}
% \inf_{\hV} \sup_{\bfSigma\in\Theta} \Ex\fnorm{\hV\hV' -\V\V' }^2
% \gtrsim1 \wedge\pth{\frac{1}{n h(\lambda)}\pth{r(s-r) + (s-r) \log
% \label{eqlbq0}
% \end{equation}
% %where $c$ is an absolute constant. % Since $\kappa> 1$ and $h$ is
%strictly increasing,
In view of the fact that\vadjust{\goodbreak} $(a \wedge b) + (c \wedge d) \geq(a \wedge
c)(b + d)$, it is sufficient to prove the following inequalities separately:
%
%
%e53 #&#
%e54 #&#
\begin{eqnarray}
\inf_{\whV} \sup_{\bolds{\Sigma}\in\Theta} \mathsf{E}\bigl\|\whV
\whV' -\mathbf {V}\mathbf{V}' \bigr\|_{\mathrm{F}}^2
&\gtrsim& r\wedge(s-r) \wedge\frac{r(s-r)}{n h(\lambda)}, \label{eqlbq01}
\\
\inf_{\whV} \sup_{\bolds{\Sigma}\in\Theta} \mathsf{E}\bigl\|\whV
\whV' -\mathbf {V}\mathbf{V}' \bigr\|_{\mathrm{F}}^2
&\gtrsim& 1 \wedge \frac{s-r}{n h(\lambda)} \log\frac{\mathrm{e}(p-r)}{s-r}. \label{eqlbq02}
\end{eqnarray}

Inequality (\ref{eqlbq01}) follows from an oracle argument: consider
the following sub-collection:
\[
\lleft\{ \mathbf{V}=\lleft[\matrix{ \mathbf{V}_1
\cr
\bolds{0}}
\rright]\dvtx\mathbf{V}_1 \in O(s,r) \rright\}.
\]
Split the data matrix according to $\mathbf{X}= [\mathbf{X}_1,
\mathbf{X}_2]$, where $\mathbf{X}_1$ consists of the first $s$ columns.
Let $\bolds{\Lambda}=
\operatorname{diag}({\lambda_1,\ldots,\lambda_s})$. Then the rows of
$\mathbf{X}_1$ and $\mathbf{X}_2$ are i.i.d. according to
$\mathcal{N}(0,\mathbf{V}_1 \bolds{\Lambda }\mathbf{V}_1' +
\mathbf{I}_s)$ and $N(0,\mathbf{I}_{p-s})$, respectively. Therefore a
sufficient statistic for estimating $\mathbf{V}$ is $\mathbf{X}_1$.
This reduces the problem to an $s$-dimensional unconstrained PCA
problem. Applying the lower bound in Theorem~\ref{thmoracle-lb} yields
(\ref{eqlbq01}).

Inequality (\ref{eqlbq02}) follows from existing results on rank-one
estimation (e.g., \cite{Birnbaum12,Vu12}). To make the argument
rigorous, we focus on the special case where
$\{{\mathbf{v}_2,\ldots,\mathbf{v}_r}\}$ are fixed to be standard
basis. Denote the following sub-collection:
%
%
%e55 #&#
\begin{equation}
\lleft\{ \mathbf{V}=\lleft[\matrix{ \mathbf{v}_1 &\bolds{0}
\cr
\bolds{0} &\mathbf{I}_{r-1} } \rright]\dvtx\mathbf{v}_1 \in
\mathbb {S}^{p-r}, \bigl|\operatorname{supp}( \mathbf{v}_1)\bigr| \leq
s-r+1 \rright\}, \label{eqq02} %\calF_0(s-r,p-r+1,1)
\end{equation}
which is well defined since $s \leq p$ by definition.
%we have assumed that $s \leq p-r-1$ in Theorem~\ref{thmlower-bound}.
%p-(r-1)$, \ie, $s\leq p$?}
Let $\mathbf{X}=[\mathbf{X}_1,\mathbf{X}_2]$, where $\mathbf{X}_1$
denotes the first $p-r+1$ columns of
$\mathbf{X}$. Then $\mathbf{X}_1$ and $\mathbf{X}_2$ consists of $n$
independent samples from
$N(0, \mathbf{I}_{p-r+1} + \lambda\mathbf{v}_1\mathbf{v}_1')$ and
$N(0,\mathbf{I}_{r-1})$, respectively.
%Restricted on the subset \eqref{eqq02}, the minimax estimation
%error of $\V$ is lower bounded by that of estimation error of $\bfv_1$
%based on $\X_1$.
Restricted on the subset (\ref{eqq02}), the minimax estimation error of
$\mathbf{V}$ is equal to the minimax estimation error of $\mathbf
{v}_1$ based on $\mathbf{X}_1$. This is equivalent to replacing the
ambient dimension $p$ by $p-r+1$ and estimating only the leading
singular vector $\mathbf {v}_1$, which is $(s-r+1)$-sparse, under the
loss
$\|\mathbf{v}_1\mathbf{v}_1'-\hat{\mathbf{v}}_1\hat{\mathbf{v}}_1'\|_{\mathrm{F}}^2$.
Applying the minimax lower bound from \cite{Vu12}, Theorem
2.1\footnote{Note that \cite{Vu12}, Theorem 2.1, for $q=0$ only
applies to the regime where $s -r \leq(p-r)/\mathrm{e}$. This does not
affect the \emph{rate} of the lower bound (\ref{eqrank1lb}) because the
minimax rate is a nondecreasing function of the sparsity $s$. Therefore
if $s - r > (p-r)/\mathrm{e}$, we can use the lower bound for $s - r =
(p-r)/\mathrm{e}$ to obtain (\ref{eqrank1lb}), since $s - r \leq p - r$
by definition.}
(see also \cite{Birnbaum12}, Theorem 2), we have
%
%
%e56 #&#
\begin{eqnarray}
\label{eqrank1lb} \inf_{\whV} \sup_{\bolds{\Sigma}\in\Theta}
\mathsf{E}\bigl\|\whV \whV' -\mathbf {V}\mathbf{V}'
\bigr\|_{\mathrm{F}}^2 &\geq& \inf_{\whV} \sup
_{\bolds{\Sigma}\in\Theta} \mathsf {E}\bigl\|\whV\whV' - \mathbf{V}
\mathbf{V}'\bigr\|^2
\nonumber\\[-8pt]\\[-8pt]
&\gtrsim& 1 \wedge\frac{(s-r)}{nh(\lambda)} \log\frac{\mathrm
{e}(p-r)}{s-r},\nonumber %% \geq\frac{c' s}{nh(\lambda)} \log\frac{\eexp p}{s},
\end{eqnarray}
completing the proof of Theorem~\ref{thmlower-bound-0}.
%where we have used $r \leq\frac{p}{2}$ implied by assumption
%Theorem~\ref{thmlower-bound-0} is now completed.
\end{pf*}

%s6.2 #&#
\subsection{\texorpdfstring{Proof of Theorem~\protect\ref{thmagg}}{Proof of Theorem 4}}
\label{seccombipf} We first state a few technical lemmas (proved in
Section~7.10 in the supplementary material~\cite{appsm}) and an
oracle upper bound (proved in Section~7.9 in the supplementary
material~\cite{appsm}), which, in view of the lower bound in
Theorem~\ref{thmoracle-lb}, gives the optimal rates of the regular PCA
problem. Some of the proofs are relegated to the supplementary material~\cite{appsm}.

%
%le2 #&#
\begin{lemma}\label{lmmquadratic}
Let $a,b,c > 0$. Then $a x^2 \leq b x + c$ implies that $x^2 \leq\frac
{b^2}{a^2} + \frac{2c}{a}$.
\end{lemma}

\begin{pf}
Since $|x - \frac{b}{2a}| \leq\frac{\sqrt{b^2+4ac}}{2a}$, we have
$x^2 \leq\frac{b^2+b^2+4ac}{2 a^2}$.
\end{pf}

%
%
%le3 #&#
\begin{lemma}\label{lmmlossi}
Let $\bolds{\Sigma}= \mathbf{I}_p + \mathbf{V}\mathbf{D}\mathbf
{V}'$. For any $\mathbf{T}\in O(p,r)$, we have
%
%
%e57 #&#
\begin{equation}
\frac{\lambda_r}{2} \bigl\|\mathbf{V}\mathbf{V}'-\mathbf{T}\mathbf
{T}'\bigr\| _{{\mathrm{F}}}^2 \leq \bigl\langle\bolds{
\Sigma}, \mathbf{V}\mathbf{V}'-\mathbf {T}\mathbf{T}'
\bigr\rangle\leq \frac{\lambda_1}{2} \bigl\|\mathbf{V}\mathbf{V}'-
\mathbf{T}\mathbf {T}' \bigr\|_{{\mathrm{F}}}^2.
\label{eqlossi}
\end{equation}
\end{lemma}

\begin{lemma}\label{lmmLM}
Let $\mathbf{K}\in\mathbb{R}^{p \times p}$ be symmetric such that
$\operatorname{Tr}(\mathbf{K}) = 0$ and $\Vert{\mathbf{K}}
\Vert_{{\mathrm{F}}}=1$. Let $\mathbf{Z}$ be $n\times p$ consisting of
independent standard normal entries. Then for any $t > 0$, we have
%
%
%e58 #&#
\begin{equation}
\mathsf{P} \biggl( \frac{1}{\sqrt{n}} \bigl| \bigl\langle\mathbf
{Z}'\mathbf{Z}, \mathbf{K} \bigr\rangle\bigr| \geq2 t +
\frac
{2t^2}{\sqrt{n}} \biggr) \leq2 \exp \bigl(-t^2 \bigr).
\label{eqLM}
\end{equation}
\end{lemma}

%
%
%le5 #&#
\begin{lemma}\label{lmmchimax}
Let ${X_1,\ldots,X_N}$ be i.i.d. such that
%
%
%e59 #&#
\begin{equation}
\label{eqchitail} \mathsf{P} \bigl( |X_1| \geq a t + b
t^2 \bigr) \leq c \exp \bigl(-t^2 \bigr),
\end{equation}
where $a,b,c > 0$. Then
% \expect{\max_{i\in[N]} |X_i|} \leq a \sqrt{\log(c N)} + \frac{a}{2
% \label{eqchimax}
%and
%
%
%e60 #&#
\begin{equation}
\mathsf{E} \max_{i\in[N]} |X_i|^2 \leq
\bigl(2 a^2 + 8 b^2 \bigr) \log(\mathrm{e} c N) + 2 b^2
\log^2 (cN). \label{eqchimaxsq}
\end{equation}
\end{lemma}

%
%
%le6 #&#
\begin{lemma}\label{lmmtref}
Let $\mathbf{E}$ be a symmetric positive definite matrix. Let $\mathbf
{F}$ be a
symmetric matrix. Then
% Let $\F$ be a symmetric matrix with eigenvalues taking values in $[-
% \begin{equation}
$| \langle\mathbf{E}, \mathbf{F}  \rangle| \leq\Vert
{\mathbf{F}} \Vert\operatorname{Tr}(\mathbf{E})$.
% \label{eqtref}
% \end{equation}
\end{lemma}

\begin{pf}
This is a special case of von Neumann's trace inequality.
\end{pf}

%
%
%le7 #&#
\begin{lemma}\label{lmmapprox}
Let $\bolds{\Theta}\in\mathcal{F}_q(s,p)$ and $k \in[p]$, where
$\mathcal{F}_q(s,p)$ is
defined in (\ref{eqrow-space}). Let $\|\bolds{\Theta}_{(i)*}\|$
denote its $i$th
largest row norm.
% Let $A \subset[p]$ denote the collection of row indices of $\Theta$
%corresponding to the $k$ largest row $\ell_2$-norm.
Then
%
%
%e61 #&#
\begin{equation}
\sum_{i > k} \|\bolds{\Theta}_{(i)*}
\|^2 \leq\frac{q}{2-q}k(s/k)^{2/q}. \label{eqapprox}
\end{equation}
\end{lemma}

\begin{pf}
By the definition of $\mathcal{F}_q(s,p)$ in (\ref{eqrow-space}), we have
\[
\sum_{i > k} \|\bolds{\Theta}_{(i)*}
\|^2 \leq s^{q/2} \sum_{i>k}
i^{-2/q} \leq s^{q/2} \int_k^{\infty}
x^{-2/q}\, \mathrm{d}x= \frac{q}{2-q}k(s/k)^{2/q}.
\]\upqed
\end{pf}

%
%
%th9 #&#
\begin{theorem}[(Oracle risk: upper bound)]\label{thmoracle-ub}
Let $p=k$ and $r \in[k]$. Let $n \geq C_0 (r + \log\lambda)$ and
$\lambda\geq C_0 \sqrt{{\log{(n)}}/{n}}$
% $nh(\lambda) \geq C_0 \log{n}$
for some sufficiently\vadjust{\goodbreak} large constant $C_0$. Let $\whV\in O(k,r)$ be
formed by the $r$ leading singular vectors of the sample covariance
matrix $\mathbf{S}$. Let $\Theta= \Theta _0(k,k,r,\lambda,\kappa)$.
Then
%
%
%e62 #&#
\begin{equation}
\sup_{\bolds{\Sigma}\in\Theta} \mathsf{E}\bigl\|\whV\whV' - \mathbf{V}
\mathbf{V}'\bigr\|_{\mathrm{F}}^2 \lesssim r \wedge(k-r)
\wedge\frac{r(k-r)}{nh(\lambda)}. \label{eqoracleub}
\end{equation}
%
%where $C$ is a constant depending only on $\kappa$.
\end{theorem}

%
% Estimating leading $r$ in $k \times k$ problem in terms of $\LL$.
% \label{lmmoracle}
%

% Let $\ntok{Z_1}{Z_N}, \ntok{W_1}{W_N} \iid N(0,1)$. Then
% \begin{equation}
% \prob{\sum_{i=1}^N Z_i W_i \geq s + s^2} \leq
% \label{eqchicross}
% \label{lmmchicross}
%

% Instead, we may need to define $\bbX_\T= \frac{1}{n}\Fnorm{\Z\T
%as follows:
% \begin{eqnarray}
% \expect{\exp\pth{\frac{\frac{1}{n}\Fnorm{\Z\T\|}^2 - r}{C}}}
%= & ~ \pth{1-\frac{2}{nC}}^{-nr/2} \exp(-r/C) \nonumber\\
%= & ~ \exp\sth{-nr/2 \log\pth{1-\frac{2}{nC}} - r/C} \\
% Therefore $\norm{\bbX_\T}_{\psi} \leq\sqrt{\frac{r}{n}}$.

\begin{pf*}{Proof of Theorem~\ref{thmagg}}
%Note that for any $a \geq0$, the estimator \eqref{eqcombi} can
%be equivalently written as the following constrained MSE estimator
% \wh\V= \argmin_{\T\in O(k,r)} \Fnorm{\S- a \T\T'},
% \label{eqcombi1}
%which gives us $\fnorm{\S- a \V\V'}^2 + a^2 \fnorm{\hV\hV' - \V
%a^2 \fnorm{\hV\hV' - \V\V'}^2 \leq\iprod{\S-a \V\V'}{\hV\hV' - \V
Before delving into the details, we give an outline of the proof as follows:
\begin{longlist}[(3)]
\item[(1)] We find a good sparse approximation of the true singular
vectors which lies in the weak-$\ell_q$ ball defined by
(\ref{eqrow-space}).

\item[(2)] We decompose the risk into a summation of three terms,
namely the \emph{approximation error}, \emph{oracle risk} and
\emph{excess risk}, the first two of which are upper bounded
in Lemma~\ref{lmmapprox} and Theorem~\ref{thmoracle-ub},
respectively.

\item[(3)] The excess risk is controlled by a careful
concentration-of-measure analysis, which forms the core of the proof.
\end{longlist}
We also remark that by (\ref{eqradius}), (\ref{eqkqs}) and condition
(\ref{eqassumption-superbigsnr}), we have
%
%
%e63 #&#
\begin{equation}
\label{eqkq-r} k_q^* \geq r.
\end{equation}
To see this, first note that $k_0^* \geq r$ by (\ref{eqradius})
directly. When $q\in(0,2)$, if $k_q^* = p$, then $k_q^*\geq r$.
Otherwise, we have
\[
k_q^*\geq s \biggl( \frac{nh(\lambda)}{r+\log(\mathrm{e}p/k_q^*)} \biggr)^{q/2} \geq s
\bigl(C_0 k_q^* \bigr)^{q/2} \geq
C_0^{q/2} s \geq r.
\]
Here the first inequality comes from (\ref{eqkqs}), the second is due
to condition (\ref{eqassumption-superbigsnr}), the third holds since
$k_q^*\geq1$ and the last holds for sufficiently large $C_0$ in view of
(\ref{eqradius}).
\begin{longlist}
\item[\textit{Step} 1: \textit{Sparse approximation}]. Fix $\mathbf
{V}\in
O(p,r)\cap\mathcal{F}_q(s,p)$. We assume that $q > 0$. Note that this
step is superfluous if $q=0$ since $\mathbf{V}$ is already sparse. Let
$k=k_q^*$ be defined in (\ref{eqkqs}). Let $\mathcal{B}(k) = \{B
\subset[p]\dvtx|B| = k\}$. Let $A \in\mathcal{B}(k)$ denote the
collection of row indices of $\mathbf{V}$ corresponding to the $k$ largest
row norm.
% Let $A \subset[p]$ denote the collection of row indices of $\Theta$
%corresponding to the $k$ largest row $\ell_2$-norm.
Put
%
%
%e64 #&#
\begin{equation}
\tilde{\bolds{\Sigma}}= \mathbf{J}_A \bolds{\Sigma}
\mathbf{J}_A + \mathbf{J}_{A^{\mathrm{c}}} = \mathbf{J}_A
\mathbf{V}\bolds{\Lambda}\mathbf{V}' \mathbf{J}_A +
\mathbf{I}_p, \label{eqtSS}
\end{equation}
where $\mathbf{J}_A$ is the diagonal matrix defined in (\ref{eqJB}).
Denote the
SVD of $\mathbf{J}_A \mathbf{V}\bolds{\Lambda}\mathbf{V}' \mathbf
{J}_A$ by $\tilde{\mathbf{V}}\tilde{\bolds{\Lambda}}\tilde
{\mathbf{V}}'$, where $\tilde{\bolds{\Lambda}}=
\operatorname{diag}({{\tilde{\lambda}_1},\ldots,{\tilde{\lambda
}_r}},{0,\ldots,0})$ and
$\tilde{\mathbf{V}}\in O(p,r) \cap\mathcal{F}_0(s,p)$, since
$\operatorname{supp}(\tilde{\mathbf{V}}) = A$. Now we claim
that $\tilde{\mathbf{V}}$ is in fact the $r$ leading singular vectors
of $\tilde{\bolds{\Sigma}}$. To
this end, note that the singular values of $\tilde{\bolds{\Sigma}}$
are $\{
{{1+\tilde{\lambda}_1},\ldots,{1+\tilde{\lambda}_r}}, 1\}$. In
view of
(\ref{eqtSS}), it is sufficient to show that the $r$th largest singular
value of $\tilde{\bolds{\Sigma}}$ is separated from one, that is,
$\sigma_{r}(\tilde{\bolds{\Sigma}}) >
1$. By Weyl's theorem (\cite{HornJohnson}, Theorem 4.3.1),
\[
\sigma_{r}(\tilde{\bolds{\Sigma}}) \geq\sigma_r(\bolds{
\Sigma}) - \|\bolds{\Sigma}- \tilde{\bolds{\Sigma}}\| \geq1+\lambda_r
- \|\bolds{\Sigma }-\tilde{\bolds{\Sigma}}\|_{\mathrm{F}}.
\]
Put $\mathbf{U}= \mathbf{J}_A \mathbf{V}$. Then
%
%
%e65 #&#
%e66 #&#
%e67 #&#
\begin{eqnarray}
\|\tilde{\bolds{\Sigma}}-\bolds{\Sigma}\|_{\mathrm{F}} &= & \bigl\|\mathbf {V}
\bolds{\Lambda}\mathbf{V}' - \mathbf{U}\bolds{\Lambda} \mathbf
{U}'\bigr\|_{\mathrm{F}}
\nonumber
\\
&\leq& \bigl\|(\mathbf{V}-\mathbf{U}) \bolds{\Lambda}\mathbf{V}'\bigr\|
_{\mathrm{F}} + \bigl\|\mathbf{U}\bolds{\Lambda}(\mathbf{V}-\mathbf{U})'
\bigr\| _{\mathrm{F}} % \nonumber\\
% %\leq& \fnorm{\V-\U} \norm{\LL} (\norm{\V} + \norm{\U}) \nonumber\\
% \leq&
\leq2 {
\lambda_1} \|\mathbf{V}- \mathbf{U}\|_{\mathrm{F}}
\nonumber
\\
&\leq& 2 \lambda_1 \sqrt{ \frac{q}{2-q} k(s/k)^{2/q}}
\label{eqcut0}
\\
&\leq& 2 \lambda_1 \sqrt{ \frac{q}{2-q} \Psi(k,p,r,n,\lambda) }
\label{eqcut1}
\\
&\leq& \frac{\lambda_r}{2}, \label{eqcut}
\end{eqnarray}
where (\ref{eqcut0}) follows from applying Lemma~\ref{lmmapprox},
(\ref{eqcut1}) follows from the choice of $k=k^*_q$ in (\ref{eqkqs}),
and (\ref{eqcut}) is implied by the assumption
(\ref{eqassumption-superbigsnr}).
%We need
%%\fnorm{\V\V'-\tV\tV'}^2 \leq\frac{q }{2 - q} k\left(r + \log
% By the choice of $k = k^*$ in \eqref{eqkqs},
Therefore
%
%
%e68 #&#
\begin{equation}
\sigma_{r}(\tilde{\bolds{\Sigma}}) \geq1 + \frac{\lambda_r}{2},
\label{eqsigmar}
\end{equation}
which implies that $\tilde{\mathbf{V}}$ indeed corresponds to the $r$
leading singular
vectors of $\tilde{\bolds{\Sigma}}$. Hence we obtain the SVD of
(\ref{eqtSS}) as
% \begin{equation}
$\tilde{\bolds{\Sigma}}= \tilde{\mathbf{V}}\tilde{\bolds{\Lambda
}}\tilde{\mathbf{V}}' + \mathbf{I}_p$.
% \label{eqtSSSVD}
% \end{equation}
Using Theorem 10 in the supplementary material~\cite{appsm}
% Theorem~\ref{thmsin-theta-cov},
we show that $\tilde{\mathbf{V}}$ provides a good sparse
approximation of $\mathbf{V}$,
%
%
%e69 #&#
\begin{equation}
\label{eqapproxerror} \bigl\|\mathbf{V}\mathbf{V}'-\tilde{\mathbf{V}}\tilde{
\mathbf{V}}'\bigr\| _{\mathrm{F}}^2 \leq \frac{2 \|\bolds{\Sigma}-\tilde{\bolds{\Sigma}}\|_{\mathrm{F}}^2}{(\sigma
_r(\tilde{\bolds{\Sigma}}) - 1)^2}
\leq\frac{32 q \kappa^2}{2-q} \Psi(k,p,r,n,\lambda),
\end{equation}
where the last inequality follows from (\ref{eqcut0}) and
(\ref{eqsigmar}). If $q = 0$, then we define $\tilde{\mathbf{V}}=
\mathbf{V}$.

\item[\textit{Step} 2: \textit{Risk decomposition}.] By definition of
    the maximizer $B^*$ in (\ref{eqBstar}), $
    \langle\mathbf{S}_{(2)},   \mathbf {V}_A \mathbf{V}_A'-\mathbf{V}_*\mathbf{V}_*'  \rangle
    \leq0$. In view of Lemma~\ref{lmmlossi}, we have
%
%
%e70 #&#
%e71 #&#
\begin{eqnarray}
 && \frac{\lambda_r}{2} \bigl\|\whV_*\whV_*'-\mathbf{V}
\mathbf{V}'\bigr\|_{\mathrm{F}}^2
\nonumber\hspace*{-30pt}
\\[-3pt]
&&\qquad\leq \bigl\langle\bolds{\Sigma}, \mathbf{V}\mathbf {V}'-
\whV_*\whV_*' \bigr\rangle
\nonumber\hspace*{-30pt}
\\
&&\qquad= \bigl\langle\bolds{\Sigma}, \mathbf{V}\mathbf {V}'-
\tilde{\mathbf{V}}\tilde{\mathbf{V}}' \bigr\rangle+ \bigl\langle
\bolds{\Sigma}, \tilde{\mathbf{V}}\tilde{\mathbf {V}}'-
\whV_A\whV_A' \bigr\rangle + \bigl\langle
\bolds{\Sigma}, \whV_A\widehat{\mathbf{V}}_A'-
\whV_* \whV_*' \bigr\rangle
\nonumber\hspace*{-30pt}
\\
&&\qquad\leq \bigl\langle\bolds{\Sigma}, \mathbf{V}\mathbf {V}'-
\tilde{\mathbf{V}}\tilde{\mathbf{V}}' \bigr\rangle+ \bigl\langle
\bolds{\Sigma}, \tilde{\mathbf{V}}\tilde{\mathbf {V}}'-
\whV_A\whV_A' \bigr\rangle + \bigl\langle
\bolds{\Sigma}-\mathbf{S}_{(2)}, \whV_A\whV_A'-
\whV_*\whV_*' \bigr\rangle
\nonumber\hspace*{-30pt}
\\
&&\qquad= \bigl\langle\bolds{\Sigma}, \mathbf{V}\mathbf {V}'-
\tilde{\mathbf{V}}\tilde{\mathbf{V}}' \bigr\rangle+ \bigl\langle
\tilde{\bolds{\Sigma}}, \tilde{\mathbf{V}}\tilde{\mathbf {V}}'-
\whV_A\whV_A' \bigr\rangle+ \bigl\langle
\bolds{\Sigma}-\mathbf{S}_{(2)}, \whV_A\whV_A'-
\whV_* \whV_*' \bigr\rangle\label{eqrisk0}\hspace*{-30pt}
\\
&&\qquad\leq\frac{\lambda_1}{2} \underbrace{\bigl\|\mathbf {V}
\mathbf{V}'-\tilde{\mathbf{V}} \tilde{\mathbf{V}}'\bigr\|
_{{\mathrm{F}}}^2}_{\mathrm{approximation\ error}}\, +\, \frac{\lambda_1}{2}
\underbrace{\bigl\|\tilde{\mathbf{V}}\tilde{\mathbf{V}}'-\widehat{
\mathbf{V}}_A\whV_A'\bigr\|_{{\mathrm{F}}}^2}_{\mathrm
{oracle\ risk}}\hspace*{-30pt}
\nonumber\\[-8pt]\label{eqriskdec} \\[-8pt]
&&\quad\qquad{}  +\, \underbrace{ \bigl\langle
\bolds{\Sigma}-\mathbf{S}_{(2)}, \whV_A\whV_A'-
\whV_*\whV_*' \bigr\rangle}_{\mathrm{excess\ risk}},\nonumber\hspace*{-30pt}
\end{eqnarray}
where (\ref{eqrisk0}) follows from that $\operatorname{supp}(\tilde
{\mathbf{V}}) = \operatorname{supp}(\whV_A) = A$, and (\ref{eqriskdec})
follows from Lemma~\ref{lmmlossi}.

Note that the expected oracle risk is upper bounded by
Theorem~\ref{thmoracle-ub} because the conditions of
Theorem~\ref{thmagg} imply those of Theorem~\ref{thmoracle-ub}. The
sparse approximation error can be upper bounded by
(\ref{eqapproxerror}). Moreover, in the exact sparse case ($q=0$), we
have $\tilde{\mathbf{V}}= \mathbf{V}$ and the approximation error is
zero.
% \nb{Will this destroy the optimal dependence on $\lambda$??}

\item[\textit{Step} 3: \textit{Excess risk}.] The hard part is to
control the
third
term (the worst-case fluctuation) in (\ref{eqriskdec}).
%, which will be shown to admit an upper bound which only depend on
%$k,p$ and $n$.
To this end, we decompose the sample covariance matrix as
\[
\mathbf{S}_{(2)}= \frac{1}{n} \mathbf{X}_{(2)}'
\mathbf{X}_{(2)}= \frac{1}{n} \bigl(\mathbf{V} \mathbf{D}
\mathbf{U}_{(2)}' + \mathbf{Z}_{(2)}'
\bigr) \bigl( \mathbf{U}_{(2)} \mathbf{D}\mathbf{V}' +
\mathbf{Z}_{(2)} \bigr).
\]
Then\vspace*{-3pt}
%
%
%e72 #&#
\begin{eqnarray}
\bolds{\Sigma}-\mathbf{S}_{(2)}&=& \mathbf{G}+ \mathbf{H},
\label{eqSdecom}
\end{eqnarray}
where\vspace*{-1pt}
%
%
%e73 #&#
%e74 #&#
\begin{eqnarray}
\mathbf{G}&\triangleq&\mathbf{V}\mathbf{D} \biggl( \frac{1}{n
}\mathbf{U}_{(2)}' \mathbf{U}_{(2)} -
\mathbf{I}_r \biggr) \mathbf {D}\mathbf{V}',
\label{eqGG}
\\
\mathbf{H}&\triangleq&\mathbf{I}_p - \frac{1}{n} \mathbf
{Z}_{(2)}'\mathbf{Z}_{(2)}- \frac{1}{n }
\mathbf{V}\mathbf{D}\mathbf{U}_{(2)}' \mathbf
{Z}_{(2)}-\frac{1}{n } \mathbf{Z}_{(2)}'
\mathbf{U}_{(2)} \mathbf{D}\mathbf{V}'. \label{eqHH}
\end{eqnarray}

We first deal the inner product with $\mathbf{G}$: write $
\langle\mathbf{G}, \whV_A\whV _A'-\whV_*\whV_*'  \rangle=
\langle\mathbf{G}, \whV_A\whV_A'-\mathbf{V}\mathbf{V}'
\rangle-  \langle\mathbf{G},
\whV_*\whV_*'-\mathbf{V}\mathbf{V}'  \rangle$. Note that
%
%
%e75 #&#
%e76 #&#
\begin{eqnarray}
\bigl\langle\mathbf{G}, \mathbf{V}\mathbf{V}' -
\whV_A\whV_A' \bigr\rangle &= & \biggl
\langle\mathbf{D} \biggl( \frac{1}{n }\mathbf{U}_{(2)}'
\mathbf{U}_{(2)} - \mathbf{I}_r \biggr) \mathbf{D},
\mathbf{V}' \bigl(\mathbf{V}\mathbf{V} ' -
\whV_A\widehat{\mathbf{V}}_A' \bigr)
\mathbf{V} \biggr\rangle
\nonumber
\\
&= & \biggl\langle\mathbf{D} \biggl( \frac{1}{n }\mathbf{U}_{(2)}'
\mathbf{U}_{(2)} - \mathbf{I}_r \biggr) \mathbf{D},
\mathbf{I}_r - \mathbf{V}' \whV_A
\whV_A'\mathbf{V} \biggr\rangle
\nonumber
\\
&\leq& \biggl\|\mathbf{D} \biggl( \frac{1}{n }\mathbf{U}_{(2)}'
\mathbf{U}_{(2)} - \mathbf{I}_r \biggr) \mathbf{D}\biggr\|
\operatorname{Tr} \bigl(\mathbf{I}_r - \mathbf{V}'
\whV_A\whV_A'\mathbf {V} \bigr)
\label{eqtref1}
\\
&\leq& \frac{\lambda_1}{2} \biggl\|\frac{1}{n }\mathbf{U}_{(2)}'
\mathbf{U}_{(2)} - \mathbf{I}_r \biggr\|\bigl\|\mathbf{V}
\mathbf{V}' - \whV_A\whV_A'
\bigr\|_{\mathrm{F}}^2, \label{eqtref2}
\end{eqnarray}
where (\ref{eqtref2}) is due to (\ref{eqlossp}) and (\ref{eqtref1}) is
a consequence of Lemma~\ref{lmmtref}, in view of the fact that
$\mathbf{I}_r - \mathbf{V}' \whV_A\whV_A'\mathbf{V}$ is symmetric
positive semi-definite while $\mathbf{D}(\frac{1}{n }\mathbf{U}_{(2)}'
\mathbf{U}_{(2)} - \mathbf{I}_r) \mathbf{D}$ is symmetric. Similarly,
we have
%
%
%e77 #&#
\begin{eqnarray}
\bigl\langle\mathbf{G}, \whV_*\whV_*'-\mathbf{V}
\mathbf{V}' \bigr\rangle&\leq& \frac{\lambda_1}{2} \biggl\|
\frac{1}{n }\mathbf{U}_{(2)}' \mathbf
{U}_{(2)} - \mathbf{I}_r \biggr\|\bigl\|\mathbf{V}
\mathbf{V}' - \widehat{\mathbf{V}}_* \whV_*'
\bigr\|_{\mathrm{F}}^2. \label{eqtref3}
\end{eqnarray}
Combining (\ref{eqtref2}) and (\ref{eqtref3}), we arrive at
%
%
%e78 #&#
\begin{eqnarray}
\bigl| \bigl\langle\mathbf{G}, \whV_A\widehat{\mathbf{V}}_A'-
\whV_*\whV_*' \bigr\rangle \bigr|&\leq& 2 \lambda_1
\biggl\|\frac{1}{n }\mathbf {U}_{(2)}'
\mathbf{U}_{(2)} - \mathbf{I}_r \biggr\|\bigl\|\whV_A
\whV_A'- \whV_*\whV_*'
\bigr\|_{\mathrm{F}}^2. \label{eqtref4}
\end{eqnarray}

Next we control the inner product with $\mathbf{H}$: recall that $A =
\operatorname{supp}(\tilde{\mathbf{V}})$ is fixed. We define a
collection of $p\times p$ symmetric
matrices indexed by $B\in\mathcal{B}(k)$ as follows:
%
%
%e79 #&#
\begin{equation}
\mathbf{K}_B \triangleq \bigl\|\whV_A
\whV_A'- \whV_B\whV_B'
\bigr\|_{{\mathrm{F}}}^{-1} \bigl( \whV_A
\whV_A'-\whV_B \whV_B'
\bigr), \label{eqKB}
\end{equation}
which has
%satisfying $\Tr(\K_B)=0$ and $\Fnorm{\K_B}=1$.
\emph{zero trace} and unit Frobenius norm.
%structure in order to obtain the desired upper bound.}
Recall that $\whV_*=\whV_{B^*}$. Then
%
%
%e80 #&#
\begin{eqnarray} \label{eqsfT}
\bigl\langle\mathbf{H}, \whV_A\whV_A'-
\whV_* \whV_*' \bigr\rangle &= & \bigl\|\whV_A
\whV_A'- \whV_*\whV_*'\bigr\|_{\mathrm{F}}
\langle\mathbf{H}, \mathbf{K}_{B^*} \rangle
\nonumber\\[-8pt]\\[-8pt]
&\leq& \bigl\|\whV_A\whV_A'-\widehat{
\mathbf{V}}_* \whV_*'\bigr\|_{\mathrm{F}} \underbrace{\max
_{B\in\mathcal{B}
(k)} \bigl| \langle\mathbf{H}, \mathbf{K}_{B}
\rangle \bigr|}_{\triangleq T}\nonumber
\end{eqnarray}
Assembling (\ref{eqSdecom}), (\ref{eqtref4}) and (\ref{eqsfT}), we can
upper bound the excess risk by
%
%
%e81 #&#
\begin{eqnarray}\label{eqsfT1}
\qquad&& \bigl\langle\bolds{\Sigma}-\mathbf{S}_{(2)}, \whV_A
\whV_A'- \whV_*\whV_*' \bigr\rangle
\nonumber
\\
&&\qquad= \bigl\langle\mathbf{G}, \whV_A\widehat{
\mathbf{V}}_A'-\whV_* \whV_*' \bigr
\rangle+ \bigl\langle\mathbf{H}, \whV_A\widehat{\mathbf{V}}_A'-
\whV_*\whV_*' \bigr\rangle
\\
&&\qquad\leq2 \lambda_1 \biggl\|\frac{1}{n }
\mathbf{U}_{(2)}' \mathbf{U}_{(2)} -
\mathbf{I}_r \biggr\|\bigl\|\whV_A \whV_A'-
\whV_*\whV_*'\bigr\| _{\mathrm{F}}^2+ T \bigl\|
\whV_A\whV_A'-\whV_* \whV_*'
\bigr\|_{\mathrm{F}}.\nonumber
\end{eqnarray}

Now we combine the risk decomposition (\ref{eqriskdec}) with the upper
bounds above to control the risk of our aggregated estimator
$\widehat{\mathbf{V}}_*$: to simplify notation, denote
\begin{eqnarray*}
%R = & \fnorm{\tV\tV'-\hV_A\hV_A'} \\
%M = & \norm{\frac{1}{n }\U_{(2)}' \U_{(2)} - \I_r}.
\delta&= & \bigl\|\whV_*
\whV_*'-\mathbf {V}\mathbf{V}'\bigr\|_{\mathrm{F}},
\qquad\Delta= \bigl\|\mathbf{V} \mathbf {V}'-\tilde{\mathbf{V}}\tilde{
\mathbf{V}}'\bigr\|_{\mathrm{F}},
\\
R &= & \bigl\|\tilde{\mathbf{V}}\tilde{\mathbf{V}}'-\whV_A
\whV_A'\bigr\|_{\mathrm{F}}, \qquad M = \biggl\|
\frac{1}{n }\mathbf{U}_{(2)}' \mathbf{U}_{(2)}
- \mathbf{I}_r\biggr\|.
\end{eqnarray*}
Assembling (\ref{eqriskdec}) and (\ref{eqsfT1}), we have
%
%
%e82 #&#
\begin{equation}
\biggl( \frac{\lambda_r }{2} - 6 \lambda_1 M \biggr)
\delta^2 \leq T \delta+ \bigl(\Delta^2 + R^2
\bigr) \biggl( \frac{\lambda_1}{2} + 6 \lambda_1 M \biggr) + T (R +
\Delta). \label{eqquad1}
\end{equation}
Introduce the event $E = \{M \leq\frac{1}{24 \kappa}\}$. By assumption
(\ref{eqassumption-lambdan}), $r \leq c'' n$ for a sufficiently small
constant $c''$. Then there exists a constant $c'>0$ only depending on
$\kappa$, such that $\frac{1}{24 \kappa} \geq2(\sqrt{r\over n} + t ) +
(\sqrt{r\over n} + t )^2$, where $t = \sqrt{\frac{\log(c' nh(\lambda))
}{n}}$. Applying Proposition~4 in the supplementary
material~\cite{appsm}
% Proposition~\ref{propwishart-bd}
yields
%
%
%e83 #&#
\begin{equation}
\mathsf{P} \bigl( E^{c} \bigr) \leq\frac{1}{c' n
h(\lambda)}.
\label
{eqwishsart-kappa}
\end{equation}
Conditioning on the event $E$ and using Lemma~\ref{lmmquadratic}, we
have
%
%
%e84 #&#
\begin{equation}
\delta^2 \leq\frac{32 T^2}{\lambda_r^2} + \frac{3\lambda_1
(\Delta^2 +
R^2) + 4 T (R + \Delta) }{\lambda_r}.
\label{eqdelta11}
\end{equation}
Recall from (\ref{eqlossp}) that the loss function is upper bounded by
$r \wedge(p-r)$. Taking expectation on both sides of (\ref{eqdelta11}),
and using (\ref{eqwishsart-kappa}) together with the Cauchy--Schwarz
inequality, we have
%to have unit coefficient before the oracle risk $\frac{R}{\lambda_r}$.
%Or is this even possible with the current method?? What we could get
%is the upper bound $\frac{\Ex R}{\lambda_r} + \frac{\Ex T^2}{
%enough.}}
%
%
%e85 #&#
%e86 #&#
\begin{eqnarray}
&& \mathsf{E} \bigl\|\whV_*\whV_*'-\mathbf {V}\mathbf{V}'
\bigr\|_{\mathrm{F}}^2
\nonumber
\\
&&\qquad\leq\frac{32 \mathsf{E} T^2 }{\lambda_r^2} + 3 \kappa \bigl(\Delta^2 +
\mathsf{E} R^2 \bigr) + \frac{ 4 \mathsf{E}[T (R + \Delta)] }{\lambda_r} + r \mathsf{P}
\bigl( E^{c} \bigr) \label{eqrisk3}
\\
&&\qquad\leq\frac{20 \mathsf{E} T^2 }{\lambda_r^2} + (3 \kappa+ 8) \bigl(\Delta ^2 +
\mathsf{E} R^2 \bigr) + \frac{r }{c' n h(\lambda)}. \label{eqrisk4}
\end{eqnarray}
In view of the oracle upper bound in Theorem~\ref{thmoracle-ub}, we
have
%
%
%e87 #&#
\begin{equation}
\mathsf{E} R^2 \leq C \biggl( r \wedge(k-r) \wedge\frac
{(k-r)r}{nh(\lambda)}
\biggr). \label{eqERq}
\end{equation}
By (\ref{eqapproxerror}), if $q>0$, the approximation is upper bounded
by
%
%
%e88 #&#
\begin{equation}
\Delta^2 \leq\frac{32 q \kappa^2}{2-q} \Psi(k,p,r,n,\lambda).
\label{eqEDeltaq}
\end{equation}
If $q = 0$, then $\Delta= 0$. To control the right-hand side of
(\ref{eqrisk4}), it boils down to upper bound $\mathsf{E} T^2 $. In the
sequel we shall prove that
%
%
%e89 #&#
\begin{equation}
\mathsf{E} T^2 \leq C (1+\lambda_1) \frac{k}{n}
\log\frac{\mathrm{e}p}{k} \label{eqETq0}
\end{equation}
for some absolutely constant $C$. Plugging (\ref{eqERq}),
(\ref{eqEDeltaq}) and (\ref{eqETq0}) into (\ref{eqrisk4}), we arrive at
%
%
%e90 #&#
%e91 #&#
\begin{eqnarray}
\qquad && \mathsf{E} \bigl\|\whV_*\whV_*'-\mathbf {V}\mathbf{V}'
\bigr\|_{\mathrm{F}}^2
\nonumber
\\
&&\qquad\leq\frac{C}{h(\lambda)} \frac{k}{n} \log\frac{\mathrm
{e}p}{k} +
\frac
{32 q \kappa^2}{2-q} \Psi(k,p,r,n,\lambda) \label{eqky1} + r \wedge
\frac{(k-r)r}{nh(\lambda)} + \frac{r}{c' n h(\lambda)}
\\
&&\qquad\leq C' \Psi(k,p,r,n,\lambda), \label{eqrisk5}
\end{eqnarray}
where the constant $C'$ only depends on $\kappa$. In the special case
of $q = 0$, the approximation error is $\Delta= 0$, which implies that
the second term in (\ref{eqky1}) is zero. Hence we have the following
stronger result:
%
%
%e92 #&#
\begin{eqnarray}\label{eqrisk6}
\mathsf{E} \bigl\|\whV_*\whV_*'-\mathbf {V}\mathbf{V}'
\bigr\|_{\mathrm{F}}^2 &\leq& \frac{C}{h(\lambda)} \frac{k}{n} \log
\frac{\mathrm
{e}p}{k} + r \wedge \frac{(k-r)r}{nh(\lambda)} + \frac{r}{c' n h(\lambda)}
\nonumber\\[-8pt]\\[-8pt]
&\leq& C' \Psi_0(s,p,r,n,\lambda),\nonumber
\end{eqnarray}
where $\Psi_0$ is defined in (\ref{eqrate-0}). Then (\ref{eqrisk5}) and
(\ref{eqrisk6}) imply the statement of the theorem for $q>0$ and $q=0$,
respectively.

To finish the proof of the theorem, it remains to establish
(\ref{eqETq0}). To this end, recall that $\mathbf{K}_B$ is symmetric
and $\operatorname{Tr}(\mathbf{K}_B)=0$. By the definitions of $T$ and
$\mathbf{H}$ in (\ref{eqsfT})~and~(\ref{eqHH}), respectively, we have
%
%
%e93 #&#
\begin{equation}
T \leq T_1 + 2 T_2, \label{eqTall}
\end{equation}
where we define
%T_1 \triangleq& \frac{1}{n} \max_{B \in\calB(k)} \left|\iprod{\Ztwo'
%%T_2 \triangleq& ~ \max_{B \in\calB(k)} \left|\iprod{\V\D\pth{
%T_2 \triangleq& ~ \frac{1}{n} \max_{B \in\calB(k)} \left|\iprod{\V
%%T_2 \triangleq& ~ \frac{1}{n} \max_{B \in\calB(k)} \left|\iprod{\V
%%T_4 \triangleq& ~ \frac{1}{n} \max_{B \in\calB(k)} \left|\iprod{
%
%
%e94 #&#
%e95 #&#
\begin{eqnarray}\quad
T_1 &\triangleq&\frac{1}{n} \max_{B \in\mathcal{B}(k)} \bigl|
\bigl\langle\mathbf{Z}_{(2)}'\mathbf{Z}_{(2)},
\mathbf {K}_{B} \bigr\rangle \bigr|\label{eqT1}
\\[-2pt]
%T_2 \triangleq& \max_{B \in\calB(k)} \left|\iprod{\V\D\pth{
T_2 &\triangleq& \frac{1}{n} \max_{B \in\mathcal{B}(k)} \bigl| \bigl\langle \mathbf{V}
\mathbf{D}\mathbf{U} _{(2)}' \mathbf{Z}_{(2)},
\mathbf {K}_{B} \bigr\rangle \bigr|= \frac{1}{n} \max
_{B \in\mathcal{B}(k)} \bigl| \bigl\langle \mathbf{Z}_{(2)}'
\mathbf{U}_{(2)} \mathbf{D}\mathbf{V}', \mathbf
{K}_{B} \bigr\rangle \bigr|\label{eqT3}. %T_2 \triangleq& \frac{1}{n} \max_{B \in
%T_4 \triangleq& \frac{1}{n} \max_{B \in\calB(k)} \left|\iprod{\Ztwo'
\end{eqnarray}
We shall prove that
%
%
%e96 #&#
%e97 #&#
\begin{eqnarray}
\mathsf{E} T_1^2 &\leq& \frac{24k}{n} \log
\frac{\mathrm{e}p}{k} + \frac{32
k^2}{n^2} \log^2 \frac{\mathrm{e}p}{k} +
\frac{62}{n}. \label{eqET1}
\\[-2pt]
\mathsf{E} T_2^2 &\leq&
\lambda_1 \biggl( \frac{40k}{n} \log\frac{\mathrm{e} p}{k} +
\frac
{24k^2}{n^2} \log^2 \frac{\mathrm{e}p}{k} + \frac{103}{n} +
\frac
{17k}{n^2} \biggr). \label{eqET3}
\end{eqnarray}
Assembling (\ref{eqTall}) with (\ref{eqET1})--(\ref{eqT3}) and
using the fact that $(a+b)^2 \leq2(a^2+b^2)$, % $(a+b+c)^2 \leq
%3(a^2+b^2+c^2)$, % $(a+b+c+d)^2 \leq4(a^2+b^2+c^2+d^2)$
we arrive at
%
%
%e98 #&#
%e99 #&#
\begin{eqnarray}
\mathsf{E} T^2 &\leq& \mathsf{E} T_1^2 + 8
\mathsf{E} T_2^2
\nonumber
\\[-2pt]
&\leq& 1500 (1+\lambda_1) \biggl( \frac{k}{n} \log
\frac{\mathrm
{e}p}{k} + \frac{k^2}{n^2} \log^2 \frac{\mathrm{e}p}{k}
\biggr)
\\[-2pt]
&\leq& 3000 (1+\lambda_1) \frac{k}{n} \log\frac{\mathrm{e}p}{k},
\label{eqETq}
\end{eqnarray}
where we used $\frac{k}{n} \log\frac{p}{k} \leq1$ implied by the
assumption (\ref{eqassumption-lambdan}).

It then remains to establish (\ref{eqET1})--(\ref{eqET3}). Note that
the collection $\{\mathbf{K}_B\dvtx B \in\mathcal{B}(k)\}$ belongs
to the
$\sigma$-algebra generated by the first sample $\mathbf{X}_{(1)}$,
which is
independent of $(\mathbf{Z}_{(2)}, \mathbf{U}_{(2)})$. By
conditioning on $\mathbf{X}_{(1)}$, we can
treat $\{\mathbf{K}_B\dvtx B \in\mathcal{B}(k)\}$ as fixed
matrices.\hspace*{230pt}\qed
\end{longlist}\noqed%
\end{pf*}

\begin{pf*}{Proof of (\ref{eqET1})}
For each fixed $B \in\mathcal{B}(k)$, $\mathbf{K}_B \indep\mathbf
{Z}_{(2)}$. Applying
Lem\-ma~\ref{lmmLM}, we have
%$\prob{\frac{1}{\sqrt{n}} |\iprod{\Z'\Z}{\K_B}| \geq2 t + \frac{2t^2}{
%
\[
\mathsf{P} \biggl( \frac{1}{\sqrt{n}} \bigl| \bigl\langle\mathbf
{Z}'\mathbf{Z}, \mathbf{K}_B \bigr\rangle\bigr| \geq2 t +
\frac
{2t^2}{\sqrt{n}} \biggr) \leq2 \exp \bigl(-t^2 \bigr).
\]
Applying Lemma~\ref{lmmchimax} with $N=|\mathcal{B}(k)| = {p\choose k}
\leq(\frac{\mathrm{e}p}{k})^k$, $a=2,b=\frac{2}{\sqrt{n}}$ and $c=2$, we have
%
%
%e100 #&#
%e101 #&#
\begin{eqnarray}
\mathsf{E} T_1^2 &\leq& \frac{1}{n} \biggl( 8
\log(2 \mathrm{e} N) + \frac
{8}{n} \bigl(\log^2 (2N) + 2 \log(2\mathrm{e}N)
\bigr) \biggr)
\\[-2pt]
&= & \frac{24}{n} \log(2 \mathrm{e} N) + \frac{8}{n^2} \log^2
(2N),
\end{eqnarray}
which implies (\ref{eqET1}).\vadjust{\goodbreak}
\end{pf*}

%inequality and $\Fnorm{\K_B}=1$, we have
%T_2^2
%Then
%= & \lambda_1 \pth{ r \frac{\Ex\|\U_{(2)}_{1 *}\|_2^2 -n}{n^2} +
%(r^2-r) \frac{\Ex|\Iprod{\U_{(2)}_{1 *}}{\U_{(2)}_{2 *}}|^2}{n^2} } \\
%= & \frac{\lambda_1 (r^2+r)}{n}.

\begin{pf*}{Proof of (\ref{eqET3})}
Fix $B \in\mathcal{B}(k)$. Since $\mathbf{U}
_{(2)} \indep\mathbf{Z}_{(2)}$, conditioned on the realization of
$\mathbf{U}_{(2)}$,
% \[
$ \langle\mathbf{V}\mathbf{D}\mathbf{U}_{(2)}' \mathbf
{Z}_{(2)}, \mathbf{K}_{B}  \rangle=  \langle\mathbf
{K}_{B} \mathbf{V}\mathbf{D}\mathbf{U} _{(2)}', \mathbf{Z}_{(2)}'
\rangle$
% \]
is distributed according to $N (0,\|\mathbf{K}_{B} \mathbf{V}\mathbf
{D}\mathbf{U} _{(2)}'\|_{\mathrm{F}}^2 )$. Therefore
\[
\bigl\langle\mathbf{V}\mathbf{D}\mathbf{U}_{(2)}' \mathbf
{Z}_{(2)}, \mathbf{K}_{B} \bigr\rangle \stackrel{(\mathrm{d})} {=}
\bigl\|\mathbf{K}_{B} \mathbf{V}\mathbf{D}\mathbf{U}_{(2)}'\bigr\| _{{\mathrm{F}}} W
\]
for some $W\sim N(0,1)$ independent of $\mathbf{U}_{(2)}$.

Using the fact that $\Vert{\mathbf{A}\mathbf{B}} \Vert_{{\mathrm{F}}}
\leq\Vert{\mathbf{A}} \Vert_{{\mathrm{F}}} \Vert{\mathbf{B}} \Vert$, we
have
%$\Fnorm{\K_{B} \V\D\U_{(2)}'} \leq\Fnorm{\K_{B}} \norm{\V} \norm{
%
\[
\bigl\|\mathbf{K}_{B} \mathbf{V}\mathbf{D}\mathbf{U}_{(2)}'\bigr\|
 _{{\mathrm{F}}} \leq \|\mathbf{K}_{B} \|_{{\mathrm{F}}}
\Vert{\mathbf {V}} \Vert\Vert{\mathbf{D}} \Vert\bigl\| \mathbf{U}_{(2)}'\bigr\|
\Vert \leq\sqrt{\lambda_1} \Vert{\mathbf{U}_{(2)}} \Vert.
\]
Consequently, $\langle\mathbf{V}\mathbf{D}\mathbf{U}_{(2)}' \mathbf
{Z}_{(2)}, \mathbf{K}_{B} \rangle$ is stochastically dominated by\break
$\sqrt{\lambda_1} \Vert{\mathbf{U}_{(2)}} \Vert|W|$. Since
$\mathbf{U}_{(2)}$ is an $n \times r$ standard Gaussian matrix, Lemma
10 in the supplementary material~\cite{appsm} yields
%we have the following large-deviation result (see, \eg, \cite[Theorem
%II.7]{Davidson01})
%
%
%e102 #&#
\begin{equation}
\mathsf{P} \bigl( \Vert{\mathbf{U}_{(2)}} \Vert\geq\sqrt{n}+
\sqrt {r}+t \bigr) \leq \exp \biggl( -\frac{t^2}{2} \biggr), \qquad t > 0.
\label{eqdz}
\end{equation}
Applying the union bound yields
\begin{eqnarray}
&& \mathsf{P} \bigl( \|\mathbf{U}_{(2)} \| |W| \geq
\sqrt {2}(\sqrt{n}+ \sqrt{r}) t + 2 t^2 \bigr)
\nonumber
\\
&&\qquad\leq\mathsf{P} \bigl( \bigl(\Vert{\mathbf{U}_{(2)}}
\Vert-\sqrt{n}-\sqrt {r}\bigr) |W| \geq2 t^2 \bigr) +\mathsf{P} \bigl( |W| \geq\sqrt{2}t \bigr)
\nonumber
\\
&&\qquad\leq\mathsf{P} \bigl( \Vert{\mathbf{U}_{(2)}} \Vert\geq
\sqrt{n}+\sqrt {r}+ \sqrt{2}t \bigr) + 2\mathsf{P} \bigl( |W| \geq\sqrt{2}t \bigr)
\nonumber
\\
&&\qquad\leq3 \exp \bigl(-t^2 \bigr),
\nonumber
\end{eqnarray}
which the last inequality follows from (\ref{eqdz}) and the
Chernoff bound $\mathsf{P} ( W \geq\sqrt{2}t  ) \leq\frac{1}{2}\exp
(-t)$. Therefore,
\[
\mathsf{P} \biggl( \frac{ \langle\mathbf{V}\mathbf{D}\mathbf
{U}_{(2)}' \mathbf{Z}_{(2)}, \mathbf{K}_{B}  \rangle}{\sqrt {\lambda_1}} \geq \sqrt{2}(\sqrt{n}+
\sqrt{r}) t + 2 t^2 \biggr) \leq3 \exp \bigl(-t^2 \bigr).
\]
Applying Lemma~\ref{lmmchimax} with $N = {p\choose k}$ yields
\[
\mathsf{E} T_2^2 \leq\frac{4\lambda_1}{n^2} \bigl(
\bigl(8+( \sqrt{n}+\sqrt{r})^2 \bigr)\log(3 eN) + 2
\log^2(3N) \bigr),
\]
which, in view of $r \leq k$, implies the desired (\ref{eqET3}).
\end{pf*}

%Note that if $\D= \sqrt{\lambda} I_r$, then the estimator
%constrained MSE estimator
% \wh\V= \argmin_{\T\in O(k,r)} \Fnorm{\S- \lambda\T\T'},
% \label{eqcombi1}
%which gives us $\fnorm{\S- \lambda\hV\hV'}^2 \leq\fnorm{\S-

%s6.3 #&#
\subsection{\texorpdfstring{Proof of Theorem~\protect\ref{thmupper-bound}}
{Proof of Theorem 7}}
We prove the theorem in three steps. First, we verify that the
``whitening'' procedure in step~3 of the reduction scheme can be
performed. Next, we investigate the signal-to-noise ratio of the
regression problem conditional on the values of $\mathbf{U}$ and
$\mathbf{Z}^0$.
Finally, we derive the desired rates by using Theorem~\ref{thmseq} and
Wedin's sin-theta theorem \cite{Wedin72}.
\begin{longlist}[($2^\circ$)]
\item[($1^\circ$)] As a first step, we verify that the ``whitening''
step is indeed possible, which requires that $\sigma_r(\mathbf{B}) > 0$.
To this end, let $J = \operatorname{supp}(\mathbf{V}^0)$. Since
$\mathbf{B}= \mathbf{U}\mathbf{D}\mathbf{V}'\mathbf{V}^0+
\mathbf{Z}^0\mathbf{V}^0$, we have
%
%
%e103 #&#
\begin{eqnarray}
\label{eqsigma-r-B} \sigma_r(\mathbf{B}) &\geq& \sigma_r
\bigl(\mathbf{U}\mathbf{D} \mathbf{V}'\mathbf{V}^0
\bigr) - \sigma_1\bigl(\mathbf{Z}^0\mathbf{V}^0
\bigr)
\nonumber
\\[-8pt]
\\[-8pt]
&\geq& \sigma_r(\mathbf{U})\sigma_r(\mathbf{D})
\sigma_r \bigl(\mathbf{V}'\mathbf{V}^0
\bigr) - \sigma_1\bigl(\mathbf{Z}^0_J\bigr).\nonumber
\end{eqnarray}

By our assumption on $\mathbf{V}^0$, condition
(\ref{eqinitial-requirement}) is satisfied with probability at least $1
- C/[nh(\lambda)]$. By Lemma 10 in the supplementary
material~\cite{appsm} and the union bound,
%
%
%e104 #&#
\begin{eqnarray}\label{eqsigma-r-U}
\sigma_r(\mathbf{U}) &\geq& \sqrt{n} \biggl(1-\sqrt
{r\over n} - \sqrt{2\log
[nh(\lambda)]\over n} \biggr),
\nonumber\\[-8pt]\\[-8pt]
\sigma_r \bigl(\mathbf{V}'\mathbf{V}^0
\bigr) &\geq& \frac{1}{2}, \qquad|J|\leq k^*_q\nonumber
\end{eqnarray}
holds with probability at least $1 - C/[nh(\lambda)]$. Note that
assumption (\ref{eqassumption-lambdan}) implies that $n\geq C_0 r$ and
that $n\geq C_0\log[nh(\lambda)]$. Thus, for sufficiently large $C_0$
in (\ref{eqassumption-lambdan}), the first inequality in
(\ref{eqsigma-r-U}) leads to $\sigma_r(\mathbf{U})\geq\frac
{2}{3}\sqrt{n}$.
% Thus we could further lower bound $\sigma_r(\U)$ by $C\sqrt{n}$.
Together with $\sigma_r(\mathbf{D}) = \sqrt{\lambda_r}$, the first
term in
(\ref{eqsigma-r-B}) is thus lower bounded by
$\frac{1}{3}\sqrt{n\lambda_r}$, and hence
%
%
%e105 #&#
\begin{equation}
\label{eqsigma-r-B-1} \sigma_r(\mathbf{B}) \geq\frac{1}{3}\sqrt{n
\lambda_r} - \sigma_1 \bigl(\mathbf{Z}_J^0
\bigr)
\end{equation}
with probability at least $1 - C/[nh(\lambda)]$.

% conditions \nb{$\log h(\lambda) / n \leq c$ and $n\geq Cr$}.

Turning to the second term in (\ref{eqsigma-r-B}), we first note that
it is upper bounded by $\max_{I\subset[p], |I| = k^*_q}
\Vert\mathbf{Z}^0_I \Vert$ conditioned on the event that $|J| \leq
k_q^*$. Note
that for any $t>0$, we have
\begin{eqnarray*}
&& \mathsf{P} \Bigl\{ \max_{I\subset[p], |I| = k^*_q} \bigl\|\mathbf
{Z}^0_I \bigr\|> \sqrt{n} + \sqrt{k^*_q} + t
\Bigr\}
\\
&&\qquad\leq\sum_{I\subset[p], |I| = k^*_q} \mathsf{P} \bigl\{ \bigl\|
\mathbf{Z}^0_I \bigr\|> \sqrt{n} + \sqrt{k^*_q}
+ t \bigr\}
\leq \pmatrix{p \cr k^*_q} \exp \bigl(-t^2/2\bigr)
\\
&&\qquad\leq \biggl(\frac{\mathrm{e}p}{k^*_q} \biggr)^{k^*_q}\exp
\bigl(-t^2/2 \bigr) = \exp \biggl(-\frac{t^2}{2} +
k^*_q\log \biggl( \frac{\mathrm{e}p}{k^*_q} \biggr) \biggr).
\end{eqnarray*}
Set $t = t^* = \sqrt{2k^*_q \log(\mathrm{e}p/k^*_q)} +
\sqrt{2\log[nh(\lambda)]}$. The rightmost side of the last display is
then bounded by $C/[nh(\lambda)]$. Thus, by (\ref{eqsigma-r-U}) and the
union bound,
%
%
%e106 #&#
\begin{equation}
\sigma_1\bigl(\mathbf{Z}^0_J\bigr) \leq
\sqrt{n} + \sqrt{k^*_q} + t^* \leq2\sqrt{n} \label{eqsigma-ZJ}
\end{equation}
with probability at least $1-C/[nh(\lambda)]$, where the last
inequality holds because the assumption (\ref{eqassumption-lambdan})
implies that $k_q^*\leq n/4$ and $t^*\leq n/2$ as long as $C_0$ is
sufficiently large.

% Here, the second inequality holds under conditions
% \nb{$k^*_q\log(ep/k^*_q)/nh(\lambda)\leq c$ and $\log n\lambda/(n

Under the assumption that $\lambda_r \geq C_0$ for some sufficiently
large $C_0 > 36 $, (\ref{eqsigma-r-B-1}) and (\ref{eqsigma-ZJ}) lead to
% \begin{equation}
% \label{eqsigma-r-B-bd}
$\sigma_r(\mathbf{B}) \geq c\sqrt{n\lambda_r} >0$ with probability
at least
$1-C/[nh(\lambda)]$.
% \end{equation}
This completes the first step in the proof.

\item[($2^\circ$)] Let ${\bar\mathbf{A}}= \tfrac{1}{\sqrt {2}}\mathbf{A}\mathbf{R}\mathbf{C}^{-1} =
\tfrac{1}{\sqrt{2}}\mathbf{D}\mathbf{U}'\mathbf{B}\mathbf
{R}\mathbf{C}^{-1} =
\tfrac{1}{\sqrt{2}}\mathbf{D}\mathbf{U}'\mathbf{L}$. Then $\bolds
{\Theta}= \mathbf{V}{\bar\mathbf{A}}$ in
(\ref{eqregression}). In the second step, we show that there exist
two constants $C_2>C_1 >0$ depending only on $\kappa$, such that
with probability at least $1-C/[nh(\lambda)]$,
%
%
%e107 #&#
\begin{equation}
\label{eqregression-snr} C_1\sqrt{n \lambda} \leq\sigma_r(\bar
\mathbf{A}) \leq\sigma_1(\bar\mathbf{A}) \leq C_2 \sqrt{n
\lambda}.
\end{equation}
To this end, note that (\ref{eqsigma-r-U}) and assumption
(\ref{eqassumption-lambdan}) imply
\[
\sigma_r({\bar\mathbf{A}}) \geq\frac{1}{\sqrt{2}}\sigma
_r(\mathbf{D}) \sigma_r(\mathbf{U}) \geq\sqrt{
\frac{n\lambda_r}{2}} \biggl(1-\sqrt {r\over n} - \sqrt
{2\log
[nh(\lambda)]\over n} \biggr) \geq C_1\sqrt{n\lambda}
\]
holds with probability at least $1-C/[nh(\lambda)]$.
% Here, the last inequality holds under conditions
% \nb{$\log h(\lambda) / n \leq c$ and $n\geq Cr$}.
Under the same assumption, Lemma~10 in the supplementary
material~\cite{appsm} implies
\[
\sigma_1({\bar\mathbf{A}}) \leq\frac{1}{\sqrt{2}}\sigma
_1(\mathbf{D}) \sigma_1(\mathbf{U}) \leq\sqrt{
\frac{n\lambda_1}{2}} \biggl(1+\sqrt {r\over n} + \sqrt
{2\log
[nh(\lambda)]\over n} \biggr) \leq C_2\sqrt{n\lambda}.
\]
Thus (\ref{eqregression-snr}) is established.

\item[($3^\circ$)] Next we show that,
%In \eqref{eqregression},
conditioned on the event that (\ref{eqregression-snr}) holds, the
signal matrix $\bolds{\Theta}$ lies in $\mathcal{F}_q(s',p)$ where
%
%
%e108 #&#
\begin{equation}
\label{eqregression-radius} s' \leq s\sigma_1^q(\bar
\mathbf{A}) \leq C s (n\lambda)^{q/2} \leq C s \bigl(nh(\lambda)
\bigr)^{q/2},
\end{equation}
where the middle inequality is due to (\ref{eqregression-snr}), the
last inequality follows from the assumption that $\lambda\geq C_0$ and
the first inequality is due to $\|\bolds{\Theta}\|_{q,w} \leq\|
\mathbf{V}\|_{q,w} \|\bar\mathbf{A}\|^q$, which is a consequence of
equation (110) in Section~7.1 of the supplementary
material~\cite{appsm}.
% \eqref{eqweaklq-A} in Section~\ref{secweaklq} in the appendix.
%To see the first inequality, note that
%(0,2)$, define $\calH_q(s,p) = \{\TT\in\reals^{p \times r}:
%verify that $\calH_{q}(s,p) \subset\calF_q(s,p) \subset
%for any $\V\in\calF_q(s,p)$ and any matrix $\A$, we have $\V\A\in

%The second inequality follows from \eqref{eqregression-snr}.

Let $k'$ be defined in (\ref{eqkqs-z}). We show that whenever
(\ref{eqregression-radius}) holds, we have
%
%
%e109 #&#
\begin{equation}
k'\leq C' k^*_q, \label{eqkpkq}
\end{equation}
where $k^*_q$ is the effective dimension defined in (\ref{eqkqs}), and
the constant $C'$ depends only on $q$. To see this, note that $k^*_q
\geq1$ by Remark~\ref{rmkeffdim}. Then\vspace*{-1.5pt} (\ref{eqkpkq})
holds trivially if $k'=1$. Next assume that $k'\geq2$. By definition,
$t_{k'-1}^{q/2} (k'-1) \leq s'$. Note that $\beta> 1$ and $t_k \geq r +
\log\frac {\mathrm{e} p}{k}$. By (\ref{eqregression-radius}), we have
$(k'-1) (r + \log\frac{\mathrm{e}p}{k'-1})^{q/2} \leq C s
(nh(\lambda))^{q/2}$. Hence $k'-1 \leq k_q^*(Cs,p,r, n,\lambda)
\leq\tau_q(C) k_q^*(s,p,\break r, n,\lambda)$, where the last inequality
follows from the third property of $k_q^*$ in Remark~\ref{rmkeffdim}.
This proves the desired~(\ref{eqkpkq}).
%when $\lambda> C_0$.

Let $E$ denote the event that both (\ref{eqsigma-r-U}) and
(\ref{eqregression-snr}) hold. Then
\begin{eqnarray*}
\mathsf{E}\bigl\|\whV\whV' -\mathbf {V}\mathbf{V}'
\bigr\|_{\mathrm{F}}^2 & =& \mathsf{E} \bigl\|\whV\whV' -
\mathbf{V}\mathbf{V}'\bigr\| _{\mathrm{F}}^2 {
\mathbf{1}_{ \{{E} \}}} + \mathsf{E} \bigl\|\whV\whV' -\mathbf{V}
\mathbf{V}'\bigr\| _{\mathrm{F}}^2 {\mathbf{1}_{ \{{E^c} \}}}
\\
& \leq&\mathsf{E}\bigl\|\whV\whV' -\mathbf {V}\mathbf{V}'
\bigr\|_{\mathrm{F}}^2 {\mathbf{1}_{ \{{E} \}}} +
\frac{Cr}{nh(\lambda)}.
\end{eqnarray*}
Here, the last inequality holds because the loss function is upper
bounded by $r$ and $\mathsf{P}(E^c) \leq C/[nh(\lambda)]$.

To further bound the first term on the rightmost hand side, we note
that $E$ is completely determined by $\mathbf{U}$ and $\mathbf{Z}^0$.
Hence, it is
nonrandom conditioned on $\mathbf{U}$ and $\mathbf{Z}^0$. Thus
\begin{eqnarray*}
\mathsf{E}\bigl\|\whV\whV' -\mathbf {V}\mathbf{V}'
\bigr\|_{\mathrm{F}}^2{\mathbf{1}_{ \{{E} \}}} & \leq &2 \mathsf{E}
\frac{\|\whT- \bolds{\Theta}\|
_{\mathrm{F}}^2}{\sigma_r^2(\bar\mathbf{A})} {\mathbf{1}_{ \{
{E} \}}} \leq\frac{C}{n\lambda}\mathsf{E}
\|\widehat{\bolds{\Theta}}- \bolds{\Theta} \|_{\mathrm{F}}^2 {
\mathbf{1}_{ \{{E} \}}}
\\
& =& \frac{C}{n\lambda} \mathsf{E} \bigl[\mathsf{E} \bigl[\| \whT- \bolds{
\Theta}\|_{\mathrm{F}}^2 {\mathbf{1}_{ \{{E} \}}} | \mathbf{U},
\mathbf{Z}^0 \bigr]{\mathbf{1}_{ \{{E} \}}} \bigr]
\\
& \leq&\frac{C}{n\lambda}\mathsf{E} \biggl[k' \biggl(r+\log
\frac{\mathrm{e}p}{k'} \biggr){\mathbf{1}_{ \{{E} \}}} \biggr]
\\
& \leq&\frac{Ck^*_q}{n \lambda} \biggl(r+\log\frac{\mathrm
{e}p}{k^*_q} \biggr).
\end{eqnarray*}
Here, the first inequality comes from
% \eqref{eqlossp} \nr{why do we need \eqref{eqlossp}?} and
Wedin's sin-theta theorem for SVD \cite{Wedin72}. The second inequality
comes from (\ref{eqregression-snr}). The second-to-last inequality
comes from Theorem~\ref{thmseq}. The last inequality holds because on
the event $E$, $k'\leq Ck^*_q$ in view of~(\ref{eqkpkq}), and $k
\mapsto k(r+\log(\mathrm{e}p/k))$ is increasing. We complete the proof by noting
that $1/\lambda\leq C/h(\lambda)$ holds since $\lambda>C_0$. The upper
bound $2 (r\wedge(p-r))$ holds in view of (\ref{eqlub}).
\end{longlist}

\begin{supplement} %[id-suppA]
\stitle{Supplement to ``Sparse PCA: Optimal rates and adaptive
estimation''}
\slink[doi]{10.1214/13-AOS1178SUPP}
\sdatatype{.pdf}
\sfilename{aos1178\_supp.pdf}
\sdescription{We provide proofs for all the remaining theoretical results in the paper. The proofs rely on
results in \cite{cover,Davidson01,Davis70,Johnstone01ss,KT59,Laurent00} and \cite{Tsybakov09}.}
\end{supplement}


%suskaldyti doi

% imsref loaded by linak, 2013-11-22 15:41:57
%
% imsref loaded by linak, 2013-11-25 18:01:49

\printaddresses

\end{document}